\documentclass[aop,MSNbibl,citesort,noMRlinks,dvips]{arximspdf}
\usepackage{graphicx}

\doi{10.1214/08-AOP398}
\volume{37}
\issue{1}
\pubyear{2009}
\firstpage{107}
\lastpage{142}

\makeatletter

\newcommand{\conf}{\mathcal{C}}
\newcommand{\e}[1]{\mathbf{#1}} % edge
\newcommand{\K}{\mathsf{K}} % Kasteleyn op
\newcommand{\Ki}{\mathsf{K}^{-1}} %inverse of K
\newcommand{\bv}{\mathbf{b}} %black vertex
\newcommand{\wv}{\mathbf{w}} %white vertex

\newcommand{\col}{x}
\newcommand{\iint}{\int\!\!\!\int}

\newcommand{\prob}[1]{\mathbb{P}[#1]} % probability
\newcommand{\probind}[2]{\mathbb{P}_{#1}[#2]} %
\newcommand{\ud}{\,d} %d of integration
\newcommand{\xx}{\mathbf{x}}

\newcommand{\tore}{\mathbb{T}^2}
\newcommand{\cercle}{\mathbb{S}^1}
\newcommand{\RR}{\mathbb{R}}
\newcommand{\ZZ}{\mathbb{Z}}
\newcommand{\CC}{\mathbb{C}}
\newcommand{\NN}{\mathbb{N}}

\newcommand{\sm}{\mathtt{t}}

\renewcommand{\Re}{\mathrm{Re}}

\newcommand{\detF}{\operatorname{Det}} %determinant Fredholm
\newcommand{\Res}{\operatorname{Res}}
\newcommand{\eqref}[1]{(\ref{#1})}

\newtheorem{lem}{Lemma}
\newtheorem{prop}{Proposition}

\newproclaim{rem}{Remark}

\makeatother

\begin{document}
\begin{frontmatter}

\title{The bead model and limit behaviors of dimer models\thanksref{T1}}
\runtitle{The bead model and limit behaviors of dimer models}
\thankstext{T1}{This work was initiated when the author was at
University Paris XI. The redaction was finished during a project at CWI,
Amsterdam, financially supported by the Netherlands Organization for
Scientific Research (NWO).}

\begin{aug}
\author[A]{\fnms{C\'edric} \snm{Boutillier}\ead[label=e1]{cedric.boutillier@upmc.fr}\corref{}}
\runauthor{C. Boutillier}
\pdfauthor{Cedric Boutillier}
\affiliation{Universit\'e Pierre et Marie Curie--Paris 6, LPMA}
\address[A]{Universit\'e Pierre et Marie Curie--Paris 6 \\
Laboratoire de Probabilit\'es et Mod\`eles Al\'eatoires \\
4 place Jussieu\\
F-75252 Paris Cedex 05\\
France\\
\printead{e1}} %adresu isvedimo komanda gale!
\end{aug}

% HISTORY:
\received{\smonth{4} \syear{2006}}
\revised{\smonth{11} \syear{2006}}

% ABSTRACT
%
\begin{abstract}
In this paper, we study the \emph{bead model}: beads are threaded on a
set of wires on the plane represented by parallel straight lines. We
add the constraint that between two consecutive beads on a wire; there
must be exactly one bead on each neighboring wire. We construct a
one-parameter family of Gibbs measures on the bead configurations that
are uniform in a certain sense. When endowed with one of these
measures, this model is shown to be a determinantal point process,
whose marginal on each wire is the \emph{sine process} (given by
eigenvalues of large hermitian random matrices). We prove then that
this process appears as a limit of any dimer model on a planar
bipartite graph when some weights degenerate.
\end{abstract}

% KEYWORDS
%
\setattribute{keyword}{AMS}{AMS 2000 subject classification.}

\end{frontmatter}

%s1 ###
\section{Introduction and presentation of the bead model}

\begin{figure}

\includegraphics{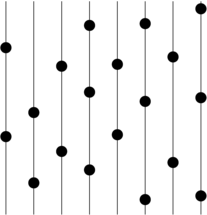}

\caption{A piece of a bead configuration.}
\label{fig:perles}
\end{figure}

We consider the collection of configurations of beads strung on an
infinite set of parallel threads lying on the plane, represented by
straight lines indexed by integers. A bead configuration on these
threads gives a configuration of points on $\ZZ\times\RR$. We impose
the following constraints on the configurations:

\begin{itemize}
\item The configuration must be locally finite: The number of beads in
each finite interval of a thread must be finite.
\item Between two consecutive beads on a thread, there must be exactly
one bead on each neighboring thread. A piece of bead configuration is
represented in Figure~\ref{fig:perles}.
\end{itemize}
Let $\Omega$ be the set of bead configurations satisfying these two conditions.

%f1 ###
%

The main goal of this paper is to construct probability measures for
our infinite system $\Omega$ that are uniform in a certain sense.

If there were only a finite number of threads of finite length, and a
fixed number of beads on each thread, then the set $\Omega$ would be a
bounded convex set of $\RR^N$, where $N$ is the total number of beads.
Therefore, the normalized Lebesgue measure on $\Omega$ would give a
uniform probability measure.
Therefore, we look for probability measures on $\Omega$ that satisfy
the two following properties:
\begin{itemize}
\item they are ergodic under the action of $\ZZ\times\RR$ by
translation,\vadjust{\goodbreak}
\item conditioned in an annular region, they induce the uniform measure
on allowed configurations inside this region.
\end{itemize}
Such a probability measure is called an \emph{ergodic Gibbs measure}.
When endowed with an ergodic Gibbs measure, the set $\Omega$ is called
a \emph{bead model}.

This model can be viewed as an interface model on $\ZZ^2$. Indeed, the
position of the beads have the same combinatorics as the square lattice
$\ZZ^2$ (rotated counterclockwise by $\pi/4$). A configuration can be
therefore encoded by a \emph{height function} $\phi\dvtx\ZZ^2
\rightarrow\RR$, where $\phi(x,y)$ is the ordinate of the bead
$(x,y)$ on its thread. The problem of existence (and uniqueness) of
Gibbs measures for this model can thus be formulated in terms of \emph
{random surfaces with a simple attractive potential}~\cite{Shef}
reflecting the hard-core interaction between beads. However, this
approach does not lead to an explicit expression for the Gibbs
measures. We adopt another point of view that will allow us to give a
closed formula for cylinder events.

% The $\sigma$-algebra $\mathcal{F}$

The $\sigma$-algebra of events for our probability measures is defined
as follows. To each bounded Borel set $B$ of $\ZZ\times\RR$ and to
each bead configuration $\omega\in\Omega$ is associated an integer
$X_B(\omega)$, equal to the number of beads in $B$. Let $\mathcal{F}$
be the smallest $\sigma$-algebra such that all the maps $X_B\dvtx
\Omega
\rightarrow\NN$ are measurable. $\mathcal{F}$ is generated by the
elementary events
\[
\{\omega\in\Omega| X_B(\omega)=n\}.
\]

%If $A$ and $B$ are disjoint Borel sets, then $X_{A\cup B}=X_A +X_B$.
%Moreover, if $(B_n)$ is a decreasing sequence of bounded Borel sets

If $\mathbb{P}$ is a Gibbs measure on $(\Omega,\mathcal{F})$, it
defines through the application $X\dvtx B\mapsto X_B$ a random process with
values in the set of boundedly finite, integer-valued measures, that
is, in other words, a \emph{random point field}.

Even without giving for the moment any explicit description for a Gibbs
measure, it is possible to estimate, at least heuristically, the
probability of some rare events. For example, one can estimate the
probability that $n$ beads lie in the same wire interval of length
$\varepsilon$. For this event to occur, then due to the geometrical
constraint imposed on configurations, there must be $n-1$ beads in a
small interval of size $\varepsilon$ on the left neighbor thread, and
$n-2$ beads in the same interval on its left, and so on, and the same
must happen on the right-hand side of the considered thread. That
is why the probability of this event must be of order $\varepsilon
^{n+2\sum_{k=1}^{n-1}(n-k)}=\varepsilon^{n^2}$. Note that this is
much smaller than the probability of having $n$ points of a \emph
{Poisson process} in such a small interval, which is of order
$\varepsilon^n$: There is thus a kind of repulsion between beads,
induced by the geometrical constraint on bead configurations.

Such a repulsion has been observed in some point processes on the real
line, especially in \emph{determinantal point processes}~\cite
{SoshDet}, for which correlation functions are expressed as
determinants of a certain kernel. Certainly, the most famous example of
such a process is the so-called \emph{sine process}, which describes
the statistics of the eigenvalues of large random hermitian matrices
with Gaussian entries (GUE ensemble~\cite{Mehta}) in the bulk of the
spectrum, and whose kernel is given by the following expression:
\[
k(x,y)=\frac{\sin(x-y)}{\pi(x-y)},
\]
in the limit when the size of the matrices goes to infinity.

We will see that the Gibbs measures we construct on $(\Omega,\mathcal
{F})$ define determinant random point fields on $\ZZ\times\RR$, for
which correlations functions are given by determinants of a kernel $J$,
whose restriction to a single thread is the sine kernel. Indeed, we
prove the following theorem:
\begin{thm}\label{thm:bead0}
For a fixed average density of beads, there exists a 1-parameter family
of ergodic Gibbs measures $(\mathbb{P}_\gamma)$ on $(\Omega,\mathcal
{F})$. When endowed with one of these measures, $(\Omega,\mathcal
{F})$ is a determinantal random point field on $\ZZ\times\RR$, with
an explicit kernel. In particular, the marginal on each thread is the
sine random point field.
\end{thm}

%f2 ###
%
\begin{figure}[b]

\includegraphics{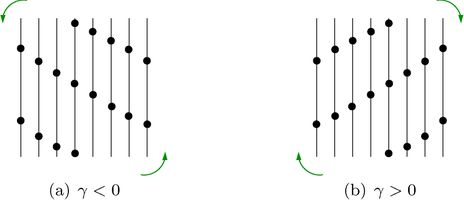}

\caption{Two typical bead configurations for different
values of $\gamma$. The parameter $\gamma$ is negative in panel
\textup{(a)},
and positive in panel \textup{(b)}. The arrows represent the effect of the
magnetic field on the beads.}
\label{fig:perles_champ}
\end{figure}

% The sine random point field process is a determinantal random point
%field on $\RR$ whose kernel is the \emph{sine kernel}
% \[
% S(x,y)=\frac{\sin(x-y)}{\pi(x-y)}
% \]
% and describes the statistics of eigenvalues in the bulk of the
%spectrum of large hermitian random matrices.
The exact expression for the kernel is given in Theorem 2. The
parameter is directly related to the average distance between a bead
and its right neighbor just below it and describes the amplitude of a
magnetic field that tends to push the beads in some direction. See
Figure \ref{fig:perles_champ}.

% take a collection of parallel wires disposed on the plane. On these
%wires, beads are threaded. We impose a constraint on configurations:

A way to construct these Gibbs measures is first to consider a
discretized version of the bead model. The set of
possible\vadjust{\goodbreak}
configurations $\Omega_\sm\subset\Omega$ is constituted by the
configurations for which the beads are located at sites of a lattice
with mesh $\sm$. We show that in this discrete setting, there exist
probability measures supported by $\Omega_\sm$, for which the
distribution of the beads is a determinantal random point field, by
exhibiting a bijection between the discretized bead model and the dimer
model on the honeycomb lattice---or equivalently random tilings of
the plane by rhombi. The measures on random tilings have the Gibbs
property and correlations have a determinantal form. Then we prove that
the sequence of discrete determinantal processes indexed by $\sm$
converges to a determinantal random point field on $\ZZ\times\RR$
when $\sm$ goes to zero.
In the second part of the paper, we prove that the bead model appears
not only as a limit of the dimer model on the honeycomb lattice, but as
the universal limit behavior of the dimer model on any bipartite
periodic planar graph.

%s2 ###
\section{A discrete version of the bead model and the dimer model on
the honeycomb lattice}

We assume for the moment that the threads are not continuous
lines, but a one-dimensional lattice with mesh size $\sm$.
The possible positions of the beads are labeled by coordinates
%
%e1 ###
%
\begin{equation}
(\col,\sm y)\in\mathbb{Z}\times\sm\mathbb{Z},
\end{equation}
$\col$ representing the thread on which the bead lies, and $y$ being
the coordinate running along the thread.

A discrete bead configuration is in bijection with a family of lattice
paths in $(\mathbb{Z}+\frac{1}{2})\times\mathbb{Z}$,
with steps going upward or to the right. A bead represents a step to
the right. Each horizontal step is then connected to its neighbor on
its right above it, by upward steps. This interpretation can be
convenient to obtain that the one-wire correlations are determinantal
and described by the sine kernel in the limit, using Lindstr\"
om--Gessel--Viennot method~\cite{GesselViennot} (or Karlin--McGregor
method~\cite{KarlinMcgregor}). See Section~\ref{sec:beadparticle}
for some complements about this interpretation.

%s2.1 ###
\subsection{Nonintersecting paths and Lindstr\"om--Gessel--Viennot
method}\label{sec:gesselviennot}

A method to study correlations between beads along a thread could have
been to look at the path interpretation of the discrete bead model in a
finite box, and then let the box grow, as suggested now.

Let$B_{n,N}$ be the finite box $\{-n+\frac{1}{2},\dots,n-\frac
{1}{2}\}\times\{N,\dots,N\}$. A~path in $B_{n,N}$ is said to be \emph
{monotonous} if its steps are going either upward or to the right.
A~family of paths $\Pi$ in $B_{n,N}$ is said to be a \emph{monotonous
$k$-path} from $(u_0,\dots,u_{k-1})$ to $(v_0,\dots,v_{k-1})$ if $\Pi
$ is a family of $k$ nonintersecting monotonous paths, such that there
exists exactly one path in $\Pi$ connecting $u_j$ to $v_j$ for every
$j\in\{0,\dots,k-1\}$.

Define $x_j=(-N,-n+j-\frac{1}{2})$ and $y_j=(N,j+\frac{1}{2})$, for
$0\leq j\leq n-1$. We endow the set of monotonous $n$-paths from
$(x_0,\dots,x_{n-1})$ to $(y_0,\dots,y_{n-1})$ with the uniform
measure. The $j$th random path will cross the vertical line $x=0$ with
a horizontal step at a random ordinate $z_j\in\{-N,\dots,N\}$. Since
the paths are not crossing one another, the points $z_j$ are distinct.

The number of monotonous paths from $x_i$ to $y_j$ is given by $
{2N+n-j+i\choose2N}$ since we need $2N$ steps upward, and $n-j+i$ to the
right. Thus, by Linstr\"om--Gessel--Viennot argument~\cite
{GesselViennot}, the number of monotonous $n$-paths from $(x_0,\dots
,x_{n-1})$ to $(y_0,\dots,y_{n-1})$
\[
\det\left[\pmatrix{2N+n-j+i\cr2N}\right]_{0\leq i,j\leq n-1}.
\]
The number of monotonous $n$-paths intersecting the vertical line $x=0$
at $z_0,\dots,z_{n-1}$ is the number of monotonous $n$-paths from
$(x_0,\dots,x_{n-1})$ to $((-\frac{1}{2}, z_0),\dots,(-\frac
{1}{2},z_{n-1}))$ times the number of monotonous $n$-paths from
$((\frac{1}{2},z_0),\dots,(\frac{1}{2},z_{n-1}))$ to
$(y_0,\dots,y_{n-1})$. As a consequence, the expression for the
probability of having crossings of the vertical axis at $z_0,\dots
,z_{n-1}$ is, again by Lindstr\"om--Gessel--Viennot's argument, a
combination of determinants
\begin{eqnarray} \label{eq:onewirecrossings}
\mathbb{P}[z_0,\dots,z_{n-1}]&=&\frac{\det[
{N+z_k+n-i-1\choose n-i-1}] \det[{N-z_k+j\choose j}
]}{\det[{2N+n-j+i\choose2N}]}\nonumber\\[-8pt]\\[-8pt]
&=&\frac{\det[
{N+z_k+i\choose i}] \det[{N-z_k+j\choose j}]}{\det
[{2N+j+i+1\choose2N}]},\nonumber
\end{eqnarray}
where the second equality is obtained from the first one by inverting
the order of the rows in the first determinant in the numerator and in
the one in the denominator.

The binomial coefficient ${N\pm z_k+i\choose i}$ is a polynomial of
degree $i$ in the variable $z_k$. Therefore, using skew-symmetry and
multilinearity of the determinant, we can replace the entries of the
$i$th row of
\[
\det\biggl[\pmatrix{N \pm z_k+i\cr i}\biggr]
\]
by any polynomial of degree $i$ evaluated at the points $z_k$, up to a
multiplicative constant. A convenient choice in this case is to use the
family of orthonormal polynomial $(\psi_i(z))_{0\leq i\leq n-1}$
obtained by Gram--Schmidt orthonormalization process from the standard
basis of polynomials $(z^i)$ with respect to the uniform measure on $\{
0,\dots,n-1\}$.
These polynomials are called the discrete Chebyshev polynomials and are
a special case of a larger family, the Hahn polynomials~(see \cite
{Nikiforov} for a reference on discrete orthogonal polynomials).

This allows us to rewrite the probability \eqref{eq:onewirecrossings} as
%
%e2 ###
%
\begin{equation}
\qquad\mathbb{P}[z_0,\dots,z_{n-1}]=\frac{1}{C_{n,N}} \det[\psi
_i(z_k)]\det[\psi_j(z_k)]=\frac{1}{C_{n,N}}\det
[k_{n,N}(z_i,z_j)],
\end{equation}
where
\[
k_{n,N}(z,z')=\sum_{j=0}^{n-1} \psi_j(z)\psi_j(z')
\]
is the self-reproducing kernel of the projector on the $n$ first
discrete Chebyshev polynomials.

Using standard techniques from orthogonal polynomials~\cite
{Nikiforov,Szego}, one can prove that the correlation functions between
intersection points for the point process $(z_i)_{0\leq i \leq n-1}$
are determinantal.

Now letting $n$ and $N$ go to infinity with the proper scaling and
making use of the asymptotic results of~\cite{BaikKriech} for
discrete orthogonal polynomials, one can prove that in the bulk [i.e.,
in a neighborhood of $(0,0)$],
the point process converges to the determinantal sine process presented above.

These techniques have been used in several papers~(see, e.g., \cite
{Joh,Joharctic} and references therein) to obtain asymptotic results on
random tilings and related models in statistical mechanics. However, in
order to get the complete determinantal behavior of the bead model, and
get an explicit expression for the kernel, we will use another method,
which will exploit a bijection between discrete bead configurations and
\emph{dimer configurations} on the honeycomb lattice $H$. But before
explaining this mapping, we recall some facts about the dimer model on~$H$.

%s2.2 ###
\subsection{The dimer model on $H$}
A \emph{dimer configuration} of the honeycomb lattice $H$ is a subset
of edges of $H$ such that every vertex is incident with exactly one
edge of this subset. A dimer configuration is also called in graph
theory a \emph{perfect matching}.
On the set of all possible dimer configurations, there is a
two-parameter family~\cite{KOS} of \emph{Gibbs probability
measures}, all elements of which satisfy the following properties:
\begin{itemize}
\item They are ergodic under the action of the lattice translations
$\ZZ\times\ZZ$.
\item If a dimer configuration is fixed in an annular region of $H$,
the dimer configurations on the inside and in the outside region are
independent, and the dimer configurations inside are uniformly distributed.
\end{itemize}

\begin{figure}

\includegraphics{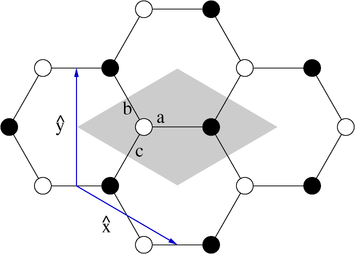}

\caption{A piece of the
honeycomb lattice $H$, with a fundamental domain (in grey), and the
system of coordinates used to label vertices. Weights $a$, $b$, $c$ are
assigned to the edges of $H$ according to their orientation.}\label
{fig:hexfunddom}
\end{figure}

Here is how these measures are defined, following~\cite
{KOS,KeLocStat}. The vertices of $H$ are colored in white and black
such that no two neighbors have the same color. Weights are assigned to
edges of $H$ according to their orientation: $a$, $b$, $c$. A
fundamental domain of $H$ is obtained by taking for example a white and
a black vertex sharing an horizontal edge. The vertices of $H$ are
named by their color (\textbf{w}hite or \textbf{b}lack) and indexed
by the coordinates of their fundamental domain $(\col,y)$. The
fundamental domain and the base vectors are represented on Figure~\ref
{fig:hexfunddom}.

%f3 ###
%

The local statistics for the Gibbs measure corresponding to these
weights have a determinantal form: The probability that edges $\e
{e}_1=(\wv_{\col_1,y_1},\bv_{\col'_1,y'_1}),\dots,\break \e{e}_k=(\wv
_{\col_k,y_k},\bv_{\col'_k,y'_k})$ belong to the random dimer
configuration are given by
%
%e3 ###
%
\begin{equation}
\prob{\e{e}_1,\dots,\e{e}_k}=\Biggl(\prod_{j=1}^k \K(\wv_{\col
_j,y_j},\bv_{\col'_j,y'_j})\Biggr)\det_{1\leq i,j\leq k} [
\Ki(\bv_{\col'_i,y'_i},\wv_{\col_j,y_j})],
\end{equation}
where $\K$ is the so-called \emph{Kasteleyn operator}. In the case of
the honeycomb lattice, it is simply the weighted adjacency matrix
restricted to the rows corresponding to white vertices, and columns
corresponding to black vertices: if $\wv$ and $\bv$ are neighbors,
then $\K(\wv,\bv)=a$, $b$ or $c$ depending on the orientation of the
edge $(\wv,\bv)$. If they are not neighbors, then $\K(\wv,\bv)=0$.
Since $\K$ is $\ZZ^2$ periodic, its inverse can be expressed using
inverse Fourier transform over the unit torus
%
%e4 ###
%
\begin{eqnarray}
\Ki(\bv_{\col',y'},\wv_{x,y})&=&\Ki\bigl(\bv_{(\col'-\col
),(y'-y)},\wv
_{0,0}\bigr)\nonumber\\[-8pt]\\[-8pt]
&=&\iint_{\tore} \frac{z^{-(y'-y)}w^{(\col'-\col
)}}{a+b/w+cz/w} \frac{\ud z}{2i\pi z} \frac{\ud w}{2i\pi w}.\nonumber
\end{eqnarray}

When none of the weights is greater than the sum of the other two, one
can show (by computing explicitly the integral above) that every type
of edge appears in the random dimer configuration with positive
probability. The measure is then said to be \emph{liquid}~\cite{KOS}.

For such weights, there is a particularly well-adapted embedding for
the hexagonal lattice, the so-called \emph{isoradial embedding}~\cite
{KeCrit} corresponding to these weights, defined as follows: The length
of a dual edge of $H$ is equal to the weight of the corresponding
primal edge, so that the dual faces of $H$ are represented as triangles
with side length $a$, $b$ and $c$. The dual faces for the usual
embedding of $H$ with regular hexagons are equilateral triangles and
the this embedding corresponds to equally distributed weights $a=b=c$.

The measure and the embedding (up to a global scaling factor) depend
only on the two ratios $a/b$ and $c/b$. For our purpose, we will chose
$a=\sm$, $b=1$, and $c=\exp(\gamma\sm)$. These values satisfy the
triangular inequalities for $\gamma\in(-1,1)$ at least for $\sm$
small enough. The distance between two successive horizontal edges is
$\sm$. We denote by $\mathbb{P}_{\gamma,\sm}$ and $\Ki_{\gamma
,\sm}$ the probability measure and the inverse Kasteleyn operator
corresponding to these particular weights.

%s2.3 ###
\subsection{Correspondence between beads and dimers}

The mapping we construct between discrete bead configurations and dimer
configurations can be described as follows: There is a bead at $(\col
,\sm y)$ if and only if the horizontal edge incident with the white
vertex in fundamental domain $(\col,y)$ is in the dimer configuration.

A way to see geometrically this correspondence is to use these
isoradial embeddings of the honeycomb lattice described above. Take an
isoradial embedding of the honeycomb lattice for weights $a=\sm$,
$b=1$, $c=e^{\gamma\sm}$ for $\gamma\in(-1,1)$ and $\sm$ small
enough. Once chosen a dimer configuration on $H$, draw a bead in the
middle of an horizontal edge if it appears in the dimer configuration
and you end with a discrete bead configuration with mesh size $\sm$.
Reciprocally, from a bead configuration, one can
reconstruct a dimer configuration by placing horizontal dimers on
edges crossing an occupied site, and completing the configuration.
This is always possible because of the intertwining of bead
positions. Moreover, the completion is unique as soon as there is at least
one bead on each wire.

%f4 ###
%
\begin{figure}

\includegraphics{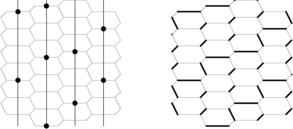}

\caption{A piece of discrete bead configuration and the
corresponding piece of dimer configuration on $H$.}
\end{figure}

% Isoradial embeddings of $H$ are parameterized, up to a global scaling
%factor, by two quantities
% $a/b,c/b$ fixing the ratios
% between the length of dual edges of each type. We take
% $a=\sm,b=1,c=e^{\gamma\sm}$, so that the distance between two
%succesive horizontal edges is $\sm$. A smaller $\sm$ corresponds to a
%finer mesh for the discrete bead model.

% The correspondance between dimer configurations and discrete bead
%configurations goes as follows:
% fix an isoradial embedding of the honeycomb
% lattice such that the middle of horizontal edges coincide with the
% possible bead positions. A bead configuration is obtained by putting
% a bead at a site each time there is an horizontal dimer at this
% position. Reciprocally, from a bead configuration one can
% reconstruct a dimer configuration by placing horizontal dimers on
% edges crossing an occupied site, and completing the configuration.
% This is always possible, because of the intertwining of bead
% positions. Moreover, the completion is unique once there is at least
% one bead on each wire.

% The parameter
% $\gamma\in(-1,1)$ parameterizes the family of isoradial embeddings
% compatible with the geometric constraints discussed above.

For a fixed $\sm$, each value of $\gamma$ corresponds to a liquid Gibbs
probability measure on dimer configurations, that can be transported
to bead configurations. The local statistics of the beads coincide with
those of
the horizontal dimers. The probability measure on bead configurations
benefits from the conditional uniform property of the dimer Gibbs measure.

This procedure defines for a given $\sm$ a family parameterized by
$\gamma$
of probability measures on discrete bead configurations that have the
conditioned uniform property.
The correlations between beads are given by determinants: the
probability of having a bead at the sites $(\col_1,\sm y_1),\dots
,(\col_k,\sm y_k)$ in the random bead configuration $\omega$ is the
probability to find the edges $(\wv_{x_1,y_1},\bv_{x_1,y_1}),\dots
,(\wv_{x_k,y_k},\bv_{x_k,y_k})$ in the random dimer configuration
%
%e5 ###
%
\begin{equation}\label{eq:probbead}
\probind{\gamma,\sm}{(\col_1,\sm y_1),\dots,(\col_k,\sm y_k) \in
\omega}=t^k \det_{1\leq i,j\leq k} \Ki_{\gamma,\sm}(\bv_{\col
_i,y_i},\wv_{\col_j,y_j}),
\end{equation}
where $\Ki_{\gamma,\sm}$ is defined by
%
%e6 ###
%
\begin{equation}\label{eq:Ki_gamma_t}
\Ki_{\gamma,\sm}(\bv_{\col,y},\wv)=\iint_{\mathbb{T}^2} \frac
{z^{-y}w^{\col}}{\sm+\frac{1}{w}(1+z e^{\gamma\sm})} \frac{\ud
z}{2i\pi z} \frac{\ud w}{2i \pi w}.
\end{equation}
%
%)} \frac{\ud u}{2i\pi} \frac{\ud w}{2i %\pi w}
% We will determine the asymptotics of this kernel for a fixed $\gamma$
%and a small $t$, and then prove the convergence of the corresponding
%determinantal process when $t$ goes to zero.

% To every compact set $B$ of
% $\mathbb{Z}\times\mathbb{R}$ is associated a random variable
% $X_K$ giving the number of beads in $K$.

% The $\sigma$-algebra of
% measurable sets is the image of that of the dimer model, and
% therefore is generated by the cylinder sets
% $\{X_{I_1}=x_1,\dots,X_{I_m}=x_m\}$ where $I_1,\dots,I_m$ are
% segments of wires.

%s3 ###
\section{Construction of explicit Gibbs measures for the continuous
bead model}\label{sec:bead_hex}

In this section, we will give an explicit description for Gibbs
measures for the bead model.
But before proving Theorem \ref{thm:bead0}, it is necessary to
investigate the behavior of the kernel defining the discrete bead
model. In other words, one has to compute the asymptotics of $\Ki
_{\gamma,\sm}(\bv_y^{\col},\wv)$ for $\sm$ small and $y$ large.

The first thing to note is that for the weights we chose, the
probability of an horizontal edge at a given white vertex (that is the
probability of a bead at a given site in the discrete model) is
%
%e7 ###
%
\begin{eqnarray}
\probind{\gamma,\sm}{\mbox{horizontal edge}}&=&\sm\Ki_{\gamma
,\sm
}(\bv_{0,0},\wv_{0,0})\nonumber\\
&=&\sm\iint_{\tore}\frac{1}{\sm+
(1/w)(1+z e^{\gamma\sm})} \frac{\ud z}{2i\pi z} \frac{\ud w}{2i
\pi w}\\
&=&\sm\sqrt{1-\gamma^2}+o(t).\nonumber
\end{eqnarray}
In order to keep a constant average density of beads, we must choose
the rescaled vertical coordinate $\xi$ equal to $\sm y \sqrt{1-\gamma^2}$.

The asymptotics of the kernel for this vertical scaling are given by
the following lemma.

% \[
% \Ki(\bv_y^{(\col)},\wv)=\iint_{\mathbb{T}^2}
%w}=\iint_{\mathbb{T}^2} \frac{u^{-\col}w^{-y}}{t+u(1+\frac{e^{
% \]
%
\begin{lem}\label{lem:noyau}
In the vertical scaling limit
$\sm\rightarrow0, \sm y \sqrt{1-\gamma^2}\rightarrow\xi$,
the coefficients $\frac{(-1)^y}{\sqrt{1-\gamma^2}}\Ki_{\gamma
,t}(\bv_{\col,y},\wv)$ converge to
%
%e8 ###
%
\begin{equation}\label{eq:kernperles}
J_{\gamma}(\col,\xi)=
\cases{
\displaystyle\int_{[-1,1]}e^{-i\xi\phi}\bigl(\gamma+i\phi\sqrt
{1-\gamma^2}\bigr)^{\col}, &\quad if $\col\geq0$, \cr
\displaystyle-\int_{\RR\setminus[-1,1]}e^{-i\xi\phi}\bigl(\gamma
+i\phi\sqrt{1-\gamma^2}\bigr)^{\col} \frac{\ud\phi}{2\pi
},&\quad if
$\col< 0$.
}
\end{equation}
In particular, when $\col=0$,
%
%e9 ###
%
\begin{equation}
J_{\gamma}(0,\xi)=\frac{1}{2\pi}\int_{[-1,1]}e^{i\xi\phi}\ud
\phi=\frac{\sin(\xi)}{\pi\xi}.
\end{equation}
\end{lem}

\begin{pf}
The entries of the inverse Kasteleyn operator are given by \eqref
{eq:Ki_gamma_t}.
%%=\iint_{\mathbb{T}^2} \frac{z^{-y}w^{x}}{a+b/w+cz} \frac{\ud z}{2i\pi
%z} \frac{\ud w}{2i \pi w}
%=\iint_{\mathbb{T}^2} \frac{z^{-y}w^{\col}}{\sm+\frac{1}{w}
%(1+z{e^{\gamma\sm}})}\frac{\ud z}{2i\pi z} \frac{\ud w}{2i
To evaluate this integral, we first perform the integration over $w$ by
the method of residues.
If $\col\geq0$, the rational fraction
% \[
% f_w(u)=\frac{u^{-\col-1}}{t+u(1+\frac{e^{\gamma t}}{w})}
% \]
%
\[
f_z(w)=\frac{w^{\col}}{ \sm w +(1+z e^{\gamma\sm})}
\]
has one pole at $w=w_0(z)=-\frac{1+z e^{\gamma\sm}}{\sm}$.
% has two poles: one at $u=0$ and the other one at $u=u_0(w)=\frac{-
% \[
% u=u_0(w)=\frac{-t}{1+\frac{e^{\gamma t}}{w}}.
% \]
% \[
% u=u_0(w)=\frac{-\sm}{1+w e^{\gamma\sm}}
% \]

%Since $u f_w(u)$ goes to zero when $u$ goes to infinity, the sum of
%the residues is zero.
By Cauchy's theorem, the integral
\[
\frac{1}{2i\pi}\int_{\mathbb{S}^1} f_z(w) \ud w
\]
is zero unless the pole $w_0(z)$ is in the unit disc, {that is,}
%
%e10 ###
%
\begin{equation}\label{eq:poleconstraint}
\Re(z) < -\frac{1+e^{2\gamma \mathtt{t}}-\mathtt{t}^2}{2e^{\gamma \mathtt{t}}}=-1+(1-\gamma^2)
{\mathtt{t}}^2 +O({\mathtt{t}}^3).
\end{equation}

Define $\theta_0=\theta_0(\gamma,\sm)=\mathrm{Arccos}(\frac
{1+e^{2\gamma\sm}-t^2}{2e^{\gamma\sm}})=\sm\sqrt{1-\gamma
^2}+O(\sm^2)$. The constraint \eqref{eq:poleconstraint} on the pole
to be inside the unit disk can be rewritten as
\[
\mathrm{arg}(z) \in(\pi-\theta_0,\pi+\theta_0).
\]
Posing $z=-e^{i\theta}=-e^{i \sm\phi\sqrt{1-\gamma^2}}$ in the
integral, we get

% \begin{align*}
% \Ki(\bv_y^{(\col)},\wv)&=
% \frac{(-1)^\col}{2\pi} t^{-\col-1}\int_{\Re(w)<-\frac{1+e^{2\gamma
%t}-t^2}{2e^{\gamma t}}} w^{-y} (1+\frac{e^{\gamma t}}{w})^
%)^\col\ud\theta\\
% &=\frac{(-1)^{y}}{2\pi}\int_{-\theta_0/t}^{\theta_0/t} e^{it y\phi}
%(\frac{e^{t(\gamma+i\phi)}-1}{t})^\col{\ud\phi}
% \end{align*}
%
\begin{eqnarray}
&&\Ki_{\gamma,\sm}(\bv_{\col,y},\wv)\nonumber\\[-8pt]\\[-8pt]
&&\qquad=
\int_{\Re(z)<-(1+e^{2\gamma\sm}-\sm^2)/(2e^{\gamma\sm})}
z^{-y}\biggl(-\frac{1+z e^{\gamma\sm}}{t}\biggr)^{\col} \frac
{\ud
z}{2\pi\sm i z} \nonumber\\
&&\qquad=(-1)^{y}\int_{-\theta_0}^{\theta_0} e^{-i y \theta}\biggl
(\frac
{e^{\gamma\sm}e^{i\theta}-1}{\sm}\biggr)^\col\frac{\ud\theta
}{2\pi\sm}\\
&&\qquad=\frac{(-1)^{y}}{2\pi}\int_{-\theta_0/\sm\sqrt{1-\gamma
^2}}^{\theta_0/\sm\sqrt{1-\gamma^2}} e^{-i\sm y\phi\sqrt{1-\gamma
^2}} \label{eq:Kiperlecpos}\nonumber\\[-8pt]\\[-8pt]
&&\quad\qquad\phantom{\frac{(-1)^{y}}{2\pi}\int_{-\theta_0/\sm
\sqrt{1-\gamma
^2}}^{\theta_0/\sm\sqrt{1-\gamma^2}}}
{}\times\biggl(\frac{e^{\sm(\gamma+i\phi\sqrt{1-\gamma
^2})}-1}{\sm
}\biggr)^\col{\sqrt{1-\gamma^2}\ud\phi}.\nonumber
\end{eqnarray}

In the vertical scaling limit $\sm\rightarrow0, t y\sqrt{1-\gamma
^2}\rightarrow\xi$, we have
\begin{eqnarray*}
\lim\frac{\theta_0}{\sm\sqrt{1-\gamma^2}} &=& 1,\\
\lim e^{-i\sm
y\phi\sqrt{1-\gamma^2}}&=&e^{-i\xi\phi},\\
\lim\frac{e^{\sm
(\gamma+i\phi\sqrt{1-\gamma^2})}-1}{\sm}&=& \gamma+i\phi\sqrt
{1-\gamma^2}.
\end{eqnarray*}

Thus, the integral above, multiplied by $\frac{(-1)^y}{\sqrt{1-\gamma
^2}}$, converges to
\[
\lim\frac{(-1)^{y}}{\sqrt{1-\gamma^2}}\Ki(\bv_{\col,y},\wv
)=\frac{1}{2\pi} \int_{[-1,1]}e^{-i\xi\phi}\bigl(\gamma+i\phi
\sqrt{1-\gamma^2}\bigr)^\col\ud\phi.
\]

When $\col< 0$, $f_z(w)$ has two poles: There is a pole at $w=0$ in
addition to that located at $w=w_0(z)=-\frac{1+z e^{\gamma\sm}}{\sm
}$. Since $w f_z(w)$ goes to zero when $z$ goes to infinity, the sum of
the residues is zero. Therefore, the integral of $f_z(w)$ on the unit
circle is not zero only if $w_0(z)$ is outside of the unit disc. It
equals in that case the opposite of the residue at $w_0(z)$. Again,
with the change of variable $z=-e^{i\theta}=-e^{i\sm\phi\sqrt
{1-\gamma^2}}$, we have
\begin{eqnarray}
\qquad\Ki_{\gamma,\sm}(\bv_{\col,y},\wv)&=&
-\int_{\Re(z)>-(1+e^{2\gamma\sm}-\sm^2)/(2e^{\gamma\sm})} z^{-y}
\biggl(-\frac{1+z e^{\gamma\sm}}{\sm}\biggr)^\col\frac{\ud
z}{2\pi\sm i z}
\\
&=&(-1)^{y+1}
\biggl(\int_{-\pi}^{-\theta_0}+\int_{\theta_0}^{\pi}\biggr)
e^{-i y \theta}\biggl(\frac{e^{\gamma\sm}e^{i\theta}-1}{\sm
}\biggr)^\col\frac{\ud\theta}{2\pi\sm}\\
&=&(-1)^{y+1}\biggl(\int_{-\pi/(\sm\sqrt{1-\gamma
^2})}^{-\theta_0/(\sm\sqrt{1-\gamma^2})} +\int_{
{\theta_0}/({\sm\sqrt{1-\gamma^2}})}^{{\pi}/{(\sm\sqrt{1-\gamma
^2})}}\biggr)e^{-i\sm y\phi\sqrt{1-\gamma^2}}\label
{eq:Kiperlecneg}\nonumber\\[-8pt]\\[-8pt]
&&
{}\times\biggl(\frac{
e^{\gamma\sm}e^{i\sm\phi\sqrt{1-\gamma^2}}-1}{\sm}\biggr)^\col
\frac{\sqrt{1-\gamma^2}\ud\phi}{2\pi}.\nonumber
\end{eqnarray}
%
% \begin{align*}
% \Ki(\bv_y^{(\col)},\wv)&=
% \frac{(-1)^\col}{2\pi} t^{-\col-1}\int_{\Re(w)>-\frac{1+e^{2\gamma
%t}-t^2}{2e^{\gamma t}}} w^{-y} (1+\frac{e^{\gamma t}}{w})^
%(1-e^{\gamma t}e^{i\theta})^\col\ud\theta\\
% &=\frac{(-1)^{y}}{2\pi}(\int_{-\pi/t}^{-\theta_0/t} +\int_{
% \end{align*}

Thus, in the scaling limit by the Lebesgue dominated convergence theorem,
\[
\lim\frac{(-1)^{y}}{\sqrt{1-\gamma^2}}\Ki_{\gamma,\sm}(\bv
_{\col,y},\wv)= \frac{-1}{2\pi}
\int_{\mathbb{R} \setminus[
-1,1]}
e^{-i\xi\phi}(\gamma+i\phi)^\col\ud\phi
\]
this completes the proof of the lemma.
\end{pf}

From the exact expressions \eqref{eq:Kiperlecpos} and \eqref
{eq:Kiperlecneg} of $\Ki(\bv^{(\col)}_y,\wv)$, one can easily check
that the entries are uniformly bounded in $y$ and $\sm$ for a given
value of $\col$, leading to the following lemma.

\begin{lem}
$\displaystyle\forall\col\in\ZZ\quad\exists M_\col>0
\quad\forall t<t_0 \quad\forall y\in\RR\quad|\Ki
_{\gamma,\sm}(\bv_{\col,y},\wv)|\leq M_\col$.
\end{lem}

These two lemmas will now be used to prove Theorem~\ref{thm:bead0},
stating that this family converges weakly to the determinantal random
point field on $\ZZ\times\RR$ with kernel~$J_{\gamma}$.
%More precisely, we will prove that for a given value of the parameter $
%$\sm$ (\ie the distributions of the horizontal edges in the dimer
%model with weights $\sm,1,e^{\gamma\sm}$) converges when $\sm$ goes to
%zero to a determinant random point field on $\ZZ\times\RR$ with the
%announced properties.

\begin{thm}\label{thm:bead1}
For each value of $\gamma\in(-1,1)$, the discrete bead model
converges weakly when $t$ goes to 0 to a determinantal random point field
on $\mathbb{Z}\times\mathbb{R}$.
The kernel of this limiting determinant random point field is $J_\gamma$:
%
%e11 ###
%
\begin{equation}
\qquad J_\gamma(\col,\xi)=
\cases{
\displaystyle\int_{[-1,1]}e^{-i\xi\phi}\bigl(\gamma
+i\phi\sqrt{1-\gamma^2}\bigr)^{\col}\frac{\ud\phi}{2\pi},
&\quad
if $\col\geq0$, \cr
\displaystyle-\int_{\RR\setminus[-1,1]}e^{-i\xi\phi}\bigl(\gamma
+i\phi\sqrt{1-\gamma^2}\bigr)^{\col}\frac{\ud\phi}{2\pi},
&\quad if $\col<0$.
}
\end{equation}
The marginal of the process on a given line is a determinantal random
point field on $\RR$ with kernel
%
%e12 ###
%
\begin{equation}
J_{\gamma}(0,\xi-\xi')=\frac{\sin(\xi-\xi')}{\pi(\xi-\xi')}.
\end{equation}
It is thus the sine random point field of the eigenvalues of
large random Hermitian matrices.
\end{thm}

\begin{pf}
Since tightness is automatic for random point fields~\cite
{pointproc}, it is
sufficient to prove the convergence of finite dimensional distributions
in order to prove the weak convergence of the family of random point
fields $(\Omega,\mathcal{F},\mathbb{P}_{\gamma,\sm})$.

Let $I_1,\dots,I_k$ be segments on wire $\col_1,\dots,\col_k,$ respectively.
It will be convenient to use multi-index notations
\[
n!=\prod_{j=1}^k n_j!,\qquad|n|=\sum_{j=1}^k
n_j,\qquad I^n=I_1^{n_1}\times\cdots\times I_k^{n_k},\qquad
z^n=z_1^{n_1}\cdots z_k^{n_k}.
\]

We will prove the convergence of the moment generating function
$G^{(I)}_{\gamma,\sm}(z_1,\dots,\break z_k)$ of the joint law of
$(X_{I_1},\dots,X_{I_k})$

%e13 ###
%
\begin{equation}
\qquad G^{(I)}_{\gamma,\sm}(z)=\mathbb{E}_{\gamma,\sm}\Biggl[\prod_{j=1}^k
(1-z_j)^{X_{I_j}}\Biggr]=\sum_{n\in\NN^k} \mathbb{E}_{\gamma,\sm
}\Biggl[\prod _{j=1}^k\frac{(X_{I_j})!}{(X_{I_j}-n_j)!}\Biggr]\frac{(-z)^{n}}{n!},
\end{equation}
where $\mathbb{E}_{\gamma.\sm}$ is the expectation with respect to
the probability $\mathbb{P}_{\gamma,\sm}$ on discrete bead configurations.
% $\mathbb{Z}\times\mathbb{R}$. To simplify, one can suppose that
% $I_1,\dots,I_k$ are segments of wires $\col_1,\dots,\col_k$
% respectively. Let
The factorial moments
%
%e14 ###
%
\begin{equation}
A^{(I)}_{\gamma,\sm}(n_1,\dots,n_k)=\mathbb{E}_{\gamma,\sm
}\Biggl[\prod _{i=1}^k \frac{X_{I_i}!}{(X_{I_i}-n_i)!}\Biggr]
\end{equation}
are quite easy to compute. They are given by the formula
%
%e15 ###
%
\begin{eqnarray}
&& A^{(I)}_{\gamma,\sm}(n_1,\dots,n_k)\hspace*{-12pt}\nonumber\\[-5pt]\\[-8pt]
&&\qquad=\mathop{\sum_{y^1_1\dots
y^1_{n_1}\in\frac{I_1}{\sm\sqrt{1-\gamma^2}}\mathrm{distinct}}}_{y^k_1\dots y^k_{n_k}\in\frac{I_k}{\sm\sqrt{1-\gamma^2}}
\mathrm{distinct}} \probind{\gamma,\sm}{\mbox{there are beads at }(\col
_1,\sm y^1_1),\dots,(\col_k,\sm y^k_{n_k})},\nonumber\hspace*{-12pt}
\end{eqnarray}
where the sum is performed over all the distinct integer $n_j$-tuples
of $\frac{I_j}{\sm\sqrt{1-\gamma^2}}$, $j=1,\dots,k$.
% \[
% p^{(I)}_{\gamma,\sm}(n_1,\dots,n_k)=\mathbb{P}_{\gamma,\sm}
%[X_{I_1}\leq n_1,\dots,X_{I_k}\leq n_k]
% \]
% is the probability of having, for all $1\leq j\leq k$, at least $n_j$
% beads in the vertical segment $I_j$.
By equation \eqref{eq:probbead}, this can be rewritten in terms of
determinants of matrices with blocks of size $n_1,\dots,n_k$
%
%e16 ###
%
\begin{equation}
\fontsize{8}{10}\selectfont{ A^{(I)}_{\gamma,\sm}(n)=\mathop{\sum_{y^1_1\dots
y^1_{n_1}\in\frac{I_1}{\sm\sqrt{1-\gamma^2}}\mathrm{distinct}}}_{
y^k_1\dots y^k_{n_k}\in\frac{I_k}{\sm\sqrt{1-\gamma^2}}\mathrm{distinct}} \sm^{|n|}\det\mathop{\left[
\begin{array}{c|c|c}
\Ki_{\gamma,\sm}(\bv_{{\col_1},y_{i1}},\wv_{{\col
_1},y_{j1}}) & \cdots&  \Ki_{\gamma,\sm}(\bv_{{\col
_k},y_{ik}},\wv_{{\col_1},y_{j1}}) \\
\hline\vdots& \ddots&  \vdots
\\
\hline\Ki_{\gamma,\sm}(\bv_{{\col_1},y_{i1}},\wv
_{{\col_k},y_{jk}}) & \cdots& \Ki_{\gamma,\sm}(\bv
_{{\col_k},y_{ik}},\wv_{{\col_k},y_{jk}})
\end{array}
\right]_{1\leq i_1,j_1 \leq n_1}}_{\hspace*{176pt} 1\leq i_k,j_k \leq
n_k}},
\end{equation}
which converges when $\sm$ goes to zero by Lemma \ref{lem:noyau} to
%e17 ###
%
\begin{eqnarray}
\quad A^{(I)}_{\gamma}(n)&=&\int_{I^n} \det\left[
\begin{array}{c|c|c}
  J_{\gamma}(\col_1-\col_1,\xi^{(1)}_{i_1}-\xi
^{(1)}_{j_1}) & \cdots&   J_{\gamma}(\col_1-\col_k,\xi
^{(1)}_{i_1}-\xi^{(k)}_{j_k})\\
\hline \vdots&  \ddots&  \vdots
\\
\hline  J_{\gamma}(\col_k-\col_1,\xi^{(k)}_{i_k}-\xi
^{(1)}_{j_1}) &\cdots&   J_{\gamma}(\col_k-\col_k,\xi
^{(k)}_{i_k}-\xi^{(k)}_{j_k})\\
\end{array}
\right]\ud^{n}\xi,\nonumber\\[-20pt]
% A^{(I)}_{\gamma}(n)=\idotsint_{I^n} \det[ J_\gamma(\hat{\xi}_j-
% ]\ud^{n}\hat\xi
\end{eqnarray}
where the integration variable $\xi$ is the $n$-tuple $(\xi
^{(1)}_1,\dots,\xi^{(1)}_{n_1},\dots,\xi^{(k)}_{n_k})$.
Since the coefficients of $\Ki_{\gamma,\sm}$ are bounded uniformly
in $\sm$ and $y$, say by $M$, then using Hadamard inequality, we get a
uniform bound on the coefficients $A^{(I)}_{\gamma,\sm}(n_1,\ldots,n_k)$
%
%e18 ###
%
\begin{equation}
\bigl|A^{(I)}_{\gamma,\sm}(n)\bigr|\leq\prod_{j=1}^k |I_j|^{n_j}
\bigl(\sqrt{|n|} M\bigr)^{|n|}.
\end{equation}
Therefore, by an argument of dominated convergence, the entire series
$Q^{(I)}_{\gamma,\sm}(z)$, $z\in\CC^k$ converges uniformly on compact
sets to
%
%e19 ###
%
\begin{equation}
Q^{(I)}_{\gamma}(z)=\sum_{n\in\NN^k} A^{(I)}_{\gamma}(n) \frac
{(-z)^{n}}{n!}
\end{equation}
which is the moment generating function for the limit distribution of
$(X_{I_1},\dots,\break X_{I_k})$.
The probability of having for all $j\in\{1,\dots,k\}$ exactly $n_j$
beads in $I_j$ is given by the following formula:
%
%e20 ###
%
\begin{equation}
\probind{\gamma}{X_{I_1}=n_1,\dots,X_{I_k}=n_k}=\frac
{(-1)^{|n|}}{n!}\cdot\frac{\partial^n}{\partial z^n}Q^{(I)}_{\gamma
}(z)\Bigm|_{z=(1,\dots,1)}.
\end{equation}

% One has to check that these quantities are finite dimensional
%distribution for some stochastic random point field.
% \begin{center}
% [Kolmogorov theorem]
% \end{center}

In particular, the probability of having no bead in a Borel set $B$ is
given by the Fredholm determinant
%
%e21 ###
%
\begin{eqnarray}
\probind{\gamma}{X_B=0}&=&\detF(\mathrm{Id}-\chi_B\Ki_{\gamma}\chi
_B)=Q^{B}_\gamma(1)\nonumber\\[-8pt]\\[-8pt]
&=&\sum_{n=0}^{\infty} \frac{(-1)^n}{n!}\int
_{B^n} \det[J_{\gamma}(\xi_i-\xi_j)]\ud^n\xi,\nonumber
\end{eqnarray}
where $\chi_B$ is the indicator function of $B$.
\end{pf}

Many quantities about the continuous bead model can be easily computed,
using the underlying dimer model. An example of such a quantity is the
average ratio between the distance between a bead and its neighbor on
the left and above it, and the distance between two successive beads on
the same thread. This average ratio is just the limit of the proportion
of $c$-edges among the nonhorizontal edges in the random dimer
configuration of $H$. It is then equal to
%
%e22 ###
%
\begin{eqnarray}
\mathbb{E}_{\gamma}[r]&=&\lim_{\sm\rightarrow0} \frac{\mathbb{P}_{\gamma
,\sm}[c\mathrm{\mbox{-}edge}]}{\mathbb{P}_{\gamma,\sm}[c\mbox{-}\mathrm{\ or\
}b\mathrm{\mbox{-}edge}]}\nonumber\\
&=&\lim_{\sm\rightarrow0} \frac{e^{\gamma\sm} \Ki
_{\gamma,\sm}(-1,-1)}{\Ki_{\gamma,\sm}(-1,0)+e^{\gamma\sm} \Ki
_{\gamma,\sm}(-1,-1)}\\
&=&\frac{\arccos{\gamma}}{\pi}.\nonumber
\end{eqnarray}
This quantity would have been difficult to obtain directly from the
description of the Gibbs measure for the point process, but has a
simple interpretation in terms of dimers.

%s4 ###
\section{Comments and interpretations of the bead model}
% %\subsection{}
%s4.1 ###
\subsection{GUE matrices and uniformly distributed intertwined points}

Afterward, it may seem not so surprising that these configurations of
``uniformly'' intertwined points are related to the determinantal sine
process and more generally to random matrices from the GUE ensemble.

Take a random (finite) matrix $H_n$ from the GUE ensemble, and define
for $k\in\{1,\dots,n\}$, $H_n^{(k)}=(H_n)_{1\leq i,j\leq k}$ to
be the submatrix of $H_n$ formed by the first $k$ rows and $k$ columns
of $H_n$. Let $\lambda_1^{(k)}\leq\lambda_k^{(k)}$ be the
eigenvalues of $H_n^{(k)}$. It follows directly from the mini-max
formulation for the eigenvalues that
%
%e23 ###
%
\begin{equation}
\forall k\in\{1,\dots,n\}, \forall j\in\{1,\dots k\}\qquad
\lambda^{(k+1)}_j\leq\lambda^{(k)}_j\leq\lambda^{(k+1)}_{j+1}.
\end{equation}
Furthermore, a result of Baryshnikov~\cite{Barysh} states that if we
condition on the eigenvalues $\lambda_1,\dots,\lambda_n$ of $H_n$,
then the eigenvalues of the submatrices are uniformly distributed over
the simplex
%
%e24 ###
%
\begin{equation}
\qquad\mathcal{S}_{(\lambda_1,\dots,\lambda_n)}=\biggl\{
\mathop{(x_j^{(k)})_{1\leq k \leq n-1}}_{\ \ \ \ \ \ \ \ \ 1\leq j\leq k} | \lambda
_j \leq x_{j}^{(n-1)}\leq\lambda_{j+1} ; x_{j}^{(k+1)}\leq
x_j^{(k)} \leq x_{j+1}^{(k+1)}\biggr\}.
\end{equation}
The bead model is somehow a bi-infinite analogue of this.

%s4.2 ###
\subsection{The bead model as an asymmetric exclusion process\label
{sec:beadparticle}}
A bead configuration can be interpreted as the history of a collection
of particles located on sites of a one-dimensional lattice $\ZZ$ and
jumping from left to right. Time is continuous and is flowing
vertically along the threads and there is a lattice site between two
successive threads. Joining every bead to the bead just above it on the
neighboring right thread, one gets an infinite collection of monotonous
paths representing the trajectories of the particles: A bead on a
thread corresponds to a jump of a particle from the site at the left of
the thread to the site on its right. Because of the geometric
constraint on beads, these paths cannot touch each other. Consequently,
the particles are submitted to an exclusion rule: a particle cannot
jump to a site if this one is already occupied by another particle.

%f5 ###
%
\begin{figure}

\includegraphics{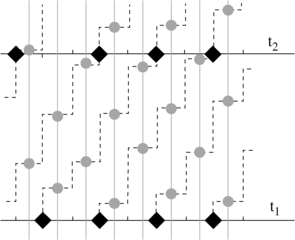}

\caption{The trajectories corresponding to the bead
configuration of Figure~\protect\ref{fig:perles} and positions of particles
(black squares) at different times $t_1$ and $t_2$.}
\end{figure}

The Gibbs measures $\mathbb{P}_{\gamma}$ on bead configurations,
viewed as families of monotonous paths constructed as above, are
probability measures on all possible evolutions of particles. The Gibbs
property and ergodicity imply that the marginal of these measures for a
fixed time (i.e., along an horizontal line) give stationary measures
for some Markovian dynamics.

The discrete bead model, as it was discussed in Section~\ref
{sec:gesselviennot}, gives a discrete version of this particle system:
In the dimer picture, a particle is represented by a $c$-edge and a
hole by a $b$-edge. Under $\mathbb{P}_{\gamma}$, the average particle
density $\rho$ is equal to the limit of the probability of a $c$-edge
and, therefore, related to $\gamma$ by the following expression
%
%e25 ###
%
\begin{equation}
\rho=\lim_{\sm\rightarrow0}\probind{\sm,\gamma}{\mbox
{$c$-edge}}=1-\frac{\arccos\gamma}{\pi}
\end{equation}
so that the density is an increasing function of $\gamma$.

The Asymmetric Simple Exclusion Process (ASEP) is also an example of
particle systems with the same constraint of exclusion. Its evolution
is Markovian, and the transition rates from an allowed configuration to
another is constant. The translation invariant stationary measures for
this model are Bernoulli probability measure, whose parameter is the density.

If the particles are not located on the vertices of the finite lattice
$\ZZ$ but on a finite annulus $\ZZ/N\ZZ$, then the number of
particles is a conserved quantity. For a fixed number of particles, the
stationary measure is uniform for ASEP. This is not the case for the
bead model with a finite number of beads\footnote{A bead model for a
finite number of threads can be constructed following the same
procedure as in the beginning of this article. We impose the number of
threads $N$ to be even to ensure that the geometric constraint on beads
makes sense. We get leading to a 1-parameter family of determinantal
random fields, whose kernel is obtained by replacing the integral in
\eqref{eq:kernperles} by a discrete sum.} the probability of a
configuration of particles depends not only on the number of particles,
but also on their positions.

The properties of the particle system coming from the bead model differ
from that of ASEP. It would be interesting to study more in details
these properties using the dimer microscopic structure, and to compare
them with that of ASEP, that are also related to random matrix
theory~\cite{Spo}.

% \subsection{The bead model and spanning trees}
%
% In \cite{Temp}, Temperley gives a bijection between dimer model on $
%correspondance has been considerably extended~\cite{KPW,KeShef}: for
%any bipartite planar graph $G$ with weighted edges, there exist a
%T-graph $G^T$ and a weight-preserving bijection between spanning trees
%of a $G^T$ and the (marked) dimer configurations of $G$. A piece of
%the T-graph corresponding to the dimer model on the honeycomb lattice
%with weights $a=\sm=0.1,b=1,c=1$, shown on figure \ref{fig:coquille},
%looks like fish scales. When $\sm$ goes to zero, the enveloppes of
%these fish scales approach a collection of superimposed cycloids.
% \begin{figure}[htp]
% \begin{center}
% \includegraphics[width=14cm]{coquille.eps}
% \caption{\footnotesize A piece of the graph on which the spanning
%tree model is in bijection with the dimer model on the honecomb
%lattice with weights $a=\sm=0.2,b=1,c=1$. The faces of this graph are
%triangles similar to the isoceles triangle with side lengths $a$, $b$
%and $c$.}\label{fig:coquille}
% \end{center}
% \end{figure}
% Using Wilson's algorithm~\cite{Wil}, one can sample a spanning tree
%on that T-graph, and thus a bead configuration using loop-erased
%random walks. It would be interesting to understand more in details
%how a random walk can sample a sine random field.

%s5 ###
\section{The bead model as a universal limit for dimer models}

Although the bead model was presented in the last section as the limit
of the dimer model on the honeycomb lattice, it turns out to be much
more general. Indeed, the bead model appears as the limit of any dimer
model on a planar $\ZZ^2$-periodic bipartite graph. We first recall
briefly some facts from the theory of the dimer model on a planar
bipartite lattice (see \cite{KOS,KO} for more details).

%s5.1 ###
\subsection{The dimer model on a bipartite planar periodic graph}
Let $G$ be a planar bipartite $\ZZ^2$-periodic graph, together with a
positive periodic weight function on the edges of $G$. We suppose that
the fundamental domain, delimited by a horizontal path $\gamma_x$ and
a vertical path $\gamma_y$, contains $n$ black vertices $\bv_1,\dots
,\bv_n$ and $n$ white vertices $\wv_1,\dots,\wv_n$, and that $G$
has at least one dimer configuration.

There is a two-parameter family of Gibbs measures on dimer
configurations of $G$ for these weights, parameterized by the two
component of an external magnetic field $B=(B_x,B_y)$~\cite{KOS}. One
can associate to this weighted graph, as in the case of the honeycomb
lattice, a \emph{Kasteleyn operator} $\K$ that will describe the
dimer model on $G$. For a given value of the magnetic field $B$, the
probability that
some edges $\e{e}_1=(\wv_1,\bv_{1}),\dots,\e{e}_k=(\wv_k,\bv_k)$
appear in the random dimer configuration is given by the following formula
%
%e26 ###
%
\begin{equation}
\probind{B}{\e{e}_1,\dots,\e{e}_k}=\Biggl(\prod_{j=1}^k \K(\wv
_j,\bv_j)\Biggr)\det_{1\leq i,j\leq k} [\Ki_B(\bv_i,\wv
_j)],
\end{equation}
where the entries of $\Ki_B$ are given by the inverse Fourier transform
%
%e27 ###
%
\begin{eqnarray}\label{eq:Kinverse}
\Ki_B(\bv^j_{x,y},\wv^i)&=&\iint_{\tore} z^{-y}w^x
[K(e^{B_x}z,e^{B_y}w)]^{-1}_{j,i} \frac{\ud z}{2 i \pi z}\frac
{\ud w}{2i\pi w}\nonumber\\[-8pt]\\[-8pt]
&=&\iint_{\tore}z^{-y}w^x\frac
{Q_{j,i}(e^{B_x}z,e^{B_y}w)}{P(e^{B_x}z,e^{B_y}w)}\frac{\ud z}{2 i \pi
z}\frac{\ud w}{2i\pi w}.\nonumber
\end{eqnarray}

The \textit{characteristic polynomial} $P(z,w)$ is the determinant of the
Fourier transform $K(z,w)$ of the periodic operator $\K$, and $Q(z,w)$
is the comatrix of $K(z,w)$.
The asymptotics of $\Ki_B,$ and thus the correlations decay depend on
the regularity of the integrand in \eqref{eq:Kinverse}, and in
particular on the presence of zeros of $P(e^{B_x}z,e^{B_y}w)$ on the
unit torus.

The \textit{spectral curve} $\{(z,w)\in\CC^2 | P(z,w)=0\}$
is a complex algebraic curve of a special kind: It is a \textit{Harnack
curve}~\cite{KO,PassRull,MikhRull}. For generic values of
$(B_x,B_y)$, $P(e^{B_x}z,e^{B_y}w)$ has zero or two roots on the unit torus,
and the phase diagram describing the behavior of the measures in
function of $B_x$ and $B_y$ is given by the \emph{am{\oe}ba} of the
\emph{spectral curve}, that is, the image of $P(z,w)=0$ by the mapping
\begin{eqnarray*}
\operatorname{Log}\dvt(\CC^*)^2&\rightarrow&\RR^2, \\
(z,w)&\mapsto&(\log|z|,\log|w|).
\end{eqnarray*}
%
% \begin{figure}[htb]
% \includegraphics[width=5cm]{amoebageneric.eps}
% \caption{\label{fig:amoebageneric}The amoeba of a Harnack curve....
%interior=complementary of convex components.}
% \end{figure}

When $(B_x,B_y)$ lies in the interior of the am{\oe}ba, the
characteristic polynomial $P(e^{B_x}z,e^{B_y}w)$ has two conjugate
roots on the unit torus and the correlations decay polynomially, and
the corresponding measure is said to be \emph{liquid} or \emph
{massless}. When $(B_x,B_y)$ lies inside a bounded connected component
of the am{\oe}ba, the correlations decay exponentially fast, and the
measure is \emph{gaseous} or \emph{massive}. When $(B_x,B_y)$ lies in
an unbounded complementary component of the am{\oe}ba, the measure is
said to be \emph{solid} and there are infinite deterministic dual
paths crossing no dimers with probability 1.

The existence of a dimer configuration on $G$ is equivalent to that of
a dimer configuration on the torus $G_1=G/\ZZ^2$. We suppose that such
a configuration on $G_1$ exists, that can be lifted to a periodic dimer
configuration $\conf_0$ of $G$.
% One can adjust the dual paths $\gamma_x$ and $\gamma_y$ delimiting
%the fundamental domain of $G$ such that they cross no dimer of $
A dimer configuration $\conf$ can be interpreted as \emph
{white-to-black unit flow}, that is, a 1-form with divergence $+1$ at
each white vertex, and divergence $-1$ at each black vertex. Therefore,
the difference $\conf-\conf_0$ is a divergence-free flow. The
horizontal (resp. vertical) \emph{slope} of a Gibbs measure $\mu$ is
the expected amount for $\mu$ of the flow $\conf-\conf_0$ across
$-\gamma_y$ (resp. $\gamma_x$). Two Gibbs measures having the same
slope are in fact equal~\cite{Shef,KOS}.
The \emph{Newton polygon} $N(P)$ of $P$, the convex hull of the
exponents of monomials of $P$, coincide with the set of all possible
slopes for a Gibbs measure on dimer configurations.
The structure of this am{\oe}ba is related to the geometry of $N(P)$
through the \emph{Ronkin function}~\cite{Mikh}
%
%e28 ###
%
\begin{equation}
R\dvtx(B_x,B_y) \mapsto\iint_{\tore} \log
(P(e^{B_x}z,e^{B_y}w)) \frac{\ud z}{2 i \pi z}\frac{\ud w}{2 i
\pi w}.
\end{equation}
In particular, the solid and gaseous phases are mapped to integer
points of $N(P)$.
% The point of $N(P)$ to which a value of the magnetic field
%$B=(B_x,B_y)$ is mapped is refered to as the \emph{slope of the Gibbs
%measure} $\mu_B$. Two Gibbs measures having the same slope are in fact
%equal~\cite{Shef,KOS}.

\begin{figure}

\includegraphics{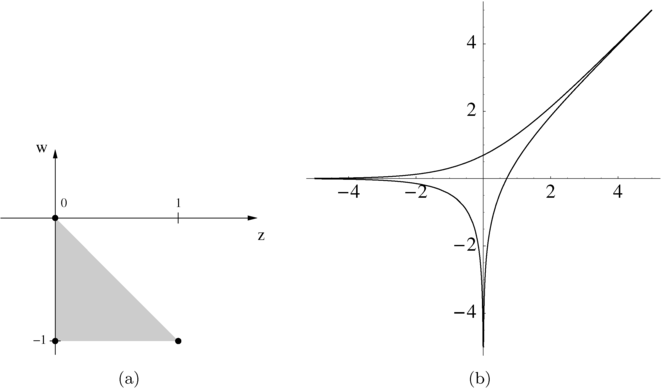}

\caption{The dimer model on the
honeycomb lattice corresponding to Figure \protect\ref{fig:hexfunddom}: (\textup{a})
the Newton polygon $N(P)$ of the characteristic polynomial
$P(z,w)=1+(1+z)/w$, and (\textup{b}) the am{\oe}ba of $P$ in the plane
$(B_x,B_y)$.}\label{fig:honeycomb}
\end{figure}

To clarify all these notions, we apply them to the particular example
of the honeycomb lattice we discussed before: Choosing weights $a$,
$b$, $c$ for edges according to their orientation (without magnetic
field) is, in fact, equivalent to choosing all weights equal to $a=\sm
$ and imposing a magnetic field equal to $B_x=\log(c/b)=\gamma\sm,
B_y=\log(a/b)=\log\sm$. The fundamental domain of the honeycomb
lattice is constituted by one black vertex and one white vertex.
Therefore, the Fourier transform $K(z,w)$ is $1\times1$ matrix:
$K(z,w)=\sm(1+1/w +z/w)$. The characteristic polynomial $P(z,w)$ is
therefore also equal to $\sm(1+1/w +z/w)$, and thus
%
%e29 ###
%
\begin{equation}
P(e^{B_x}z,e^{B_y}w)=\sm+\frac{1}{w}(1+ze^{\gamma\sm}).
\end{equation}
The cofactor $Q(z,w)$ of $K(z,w)$ is by convention equal to 1.
The Newton polygon and the am{\oe}ba for this model are represented on
Figure~\ref{fig:honeycomb}.

%f6 ###
%

The bead model is obtained from the dimer model on the honeycomb
lattice when the magnetic field goes inside one of the outgrowths of
the am{\oe}ba, in the thin region separating two solid phases.
%representing all the possible \textit{slopes} of a Gibbs measure.
In the general dimer model, the bead model will also appear near the
frontier between the liquid phase and solid phases. But before
explaining how to find the bead model in this setting, we need some
more information about the local geometry of the am{\oe}ba, in
particular about its unbounded outgrowths: the \emph{tentacles}.

\subsection{Tentacles of the am{\oe}ba}

Consider a particular side of the Newton polygon $N(P)$. Changing the
generators of the $\ZZ^2$ lattice acting on $G$ by translation induces
a linear transformation of $N(P)$. A change of basis of $\ZZ^2$ is
encoded by an element $M$ of $\mathrm{SL}_2(\ZZ)$. The linear
transformation acting on $N(P)$ is $(M^{-1})^T$. After
possibly such an operation, we can assume that this side is
horizontal, and that the polygon lies above it. Recall that the
Newton polygon represents the possible slopes for a Gibbs measure on
the dimer model on $G$. When the slope of a Gibbs measure is a lattice point
of the boundary of $N(P)$, the system is a solid phase.

We want to investigate the geometry of the phase diagram for values of the
magnetic field inducing measures with a slope close to the particular
side of $N(P)$ we chose. In particular, we seek for the shape of the
boundary of the am{\oe}ba in a neighborhood of the frontier between
the liquid phase
and the different solid phases, corresponding to the points of the
particular side of the polygon.

To get a measure with a slope close to that side of $N(P)$, we apply
to the system a magnetic field oriented essentially downward
$(B_x,B_y)=(c,-R)$, with $R$ very large. To remain close to the
notations used in the
previous section, we introduce the small parameter $\sm=e^{-R}$.

When $\sm$ is small, the leading terms in the characteristic
polynomial $P(e^c z, \sm w)$ are those with the smallest power in
$w$, say $\delta_0$.
%
%e30 ###
%
\begin{equation}
P(e^c z,\sm w) = (\sm w)^{\delta_0} \Biggl( \sum_{\gamma}
a_{\gamma\delta_0} (e^cz)^\gamma+ O(\sm) \Biggr).
\end{equation}
By a suitable choice of the
origin of the Newton polygon (deforming the paths $\gamma_x$ and/or
$\gamma_y$ delimiting the fundamental domain of $G$), one can assume
that $\delta_0=0$ and
that all the roots of $P_0(X)=\sum_{\gamma} a_{\gamma0}X^\gamma$
are positive~\cite{KOS}.

If $e^c$ is not a root of $P_0(X)=\sum_{\gamma} a_{\gamma
0}X^\gamma$, then for $\sm$ small enough, $P(e^c z,\break e^{-R} w)$ has no
roots on the unit torus. In this case, the magnetic field
$(B_x,B_y)=(c,-R)$ belongs to an unbounded component of the
am{\oe}ba. The corresponding measure $\mu(B_x,B_y)$ is solid. On
the contrary, if $e^c$ is a root of this polynomial, then for every
$R$ large enough, the polynomial has two complex conjugated roots on
the unit torus: We are in the liquid phase. The am{\oe}ba defining
the liquid phase has therefore \emph{tentacles} going to infinity
with asymptotes given by the straight lines $x=c$.

For generic weights, the
asymptotes are all distinct, and there is one asymptote for
each segment between two lattice points on the side of $N(P)$.
Moreover, one can give an asymptotic expansion for the equation of
the boundary of the am{\oe}ba. Since the boundary of the am{\oe}ba
is the image of the real locus of the curve, it is given by the
equation
%
%e31 ###
%
\begin{equation}\label{eq:reallocus}
P(\pm e^{B_x},\pm e^{B_y})=0.
\end{equation}
In a neighborhood of $(B_x,B_y)=(c,-\infty)$, the solution of \eqref
{eq:reallocus} for $B_x$ admits an
asymptotic expansion in $\sm=e^{-B_y}$: $B_x=c+c_1\sm+O(\sm^2)$.
Since $P(e^c,0)=0$, we have
%
%e32 ###
%
\begin{eqnarray}\label{eq:asympP}
P(e^{B_x},\pm e^{B_y})&=&P\bigl(e^{c+c_1\sm+O(\sm^2)},\pm\sm\bigr)\nonumber\\[-8pt]\\[-8pt]
&=&\bigl(c_1
e^c\partial_1 P(e^c,0)\pm\partial_2 P(e^c,0)\bigr)\sm+O(\sm^2).\nonumber
\end{eqnarray}
Therefore, the coefficient $c_1$ in the expansion is defined by
%
%e33 ###
%
\begin{equation}
c_1=\pm\frac{\partial_2 P(e^c,0)}{e^c \partial_1 P(e^c,0)}
\end{equation}
and the two curves $B_x=c\pm e^{-c}\frac{\partial_2 P(e^c,0)}{
\partial_1 P(e^c,0)}e^{B_y}$ define the two\vspace*{2pt} asymptotic branches of the
boundary of the am{\oe}ba in the neighborhood of
$(c,-\infty)$. Define $\beta$ as
%
%e34 ###
%
\begin{equation}
\beta= - e^{-c}\frac{\partial_2 P(e^c,0)}{\partial_1 P(e^c,0)}.
\end{equation}
For any $\gamma\in(-1,1)$, the curve
%
%e35 ###
%
\begin{equation}
B_x=c+\gamma\beta e^{B_y}
\end{equation}
lies inside the amoeba for $B_y$ negative enough.
\eject

%f7 ###
%
\begin{figure}

\includegraphics{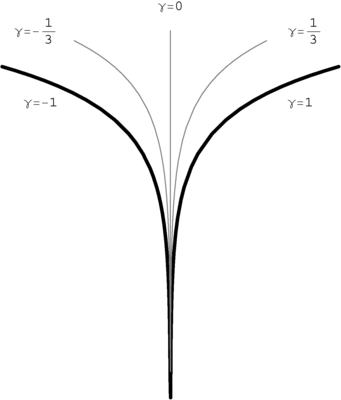}

\caption{A tentacle (in thick lines) and curves
described above, for different values of $\gamma$.}
\label{fig:tentacle}
\end{figure}

%s5.3 ###
\subsection{Deep inside a tentacle}

%s5.3.1 ###
\subsubsection{Analytic results about the roots of $P$}
Let us fix $c$ to be equal, as before to the logarithm of one of the
roots of
$P_0$. For a fixed $z$, the polynomial $P(e^{c+\gamma\beta\sm} z,
W)$ has $d$ roots $W_0(z),\dots,W_{d-1}(z)$. Since $e^c$ is a root
of $P_0$, one of these $W_j(z)$, say it is $W_0$, equals 0 when
$z=1$. If all the roots of $P_0$ are distinct, $W_0(z)$ is the only
zero having this property.
% : if there were another root $W_j(z),j\neq
% 0$ such that $W_j(1)=0$ then for $z$ close to 1, we would have more
% than two roots with the same module, what would be in contradiction
% with the 2-to-1 property of the mapping from the spectral curve to
% its am{\oe}ba.
The 2-to-1 property of the map from the spectral curve to its am{\oe
}ba shows that $W_0(z)$ does not
equal zero for $z\neq1$.
Therefore, by compactness of $\cercle$ there exists an
$\epsilon>0,$ such that
%
%e36 ###
%
\begin{equation}\label{eq:P_root_Wj}
\forall j\in\{1,\dots,d-1\}, \forall z\in\cercle\qquad |W_j(z)|\geq
\epsilon.
\end{equation}

Differentiability of the roots with respect to the coefficients and
symmetry by complex conjugation imply that there exists $\theta(\sm)=
\theta_0 \sm+O(\sm^2)$, with $\theta_0>0,$ such that
%
%e37 ###
%
\begin{equation}
|W_0(z)|\leq\sm\quad\Leftrightarrow\quad z\in
\bigl[e^{-i\theta(\sm)},e^{i\theta(\sm)}\bigr].
\end{equation}
In fact, an expansion of $P$ similar to \eqref{eq:asympP} shows that
when $\arg(z)=O(\sm)$
%
%e38 ###
%
\begin{equation}
W_0(z)=\gamma\sm+\frac{i \arg(z)}{\beta}+O(\sm^2),
\end{equation}
and, therefore, $\theta(\sm)=\sm\beta\sqrt{1-\gamma^2}+O(\sm^2)$.

%s5.3.2 ###
\subsubsection{Asymptotics of the inverse Kasteleyn operator inside a tentacle}

The Newton polygon $N(P)$ can also be obtained indirectly from the
\emph{max flow--min cut} theorem as the intersections of half-planes
made of points satisfying linear constraints in order to be a slope of
a measure on dimer coverings~\cite{KOS}.

% Recall that the Newton polygon $N(P)$ is an intersection of
%half-planes, where all the inequalities \eqref{eq:slopeconstr} are
%satisfied by the average slope of the Gibbs measure.

The side of the Newton polygon we are looking at is a segment of the
straight line delimiting one of these half-planes. When the average
slope of the Gibbs measure lies on this line, one of the inequalities
defining $N(P)$ becomes an equality, implying the existence of \emph
{frozen dual paths}, with a direction perpendicular to the side of
$N(P)$, that are not crossed by any dimer with probability 1.
%at there exists two translates of a face between which the difference
%of height equals the capacity of a dual path $\Gamma$ joining them.
%All the translated of this path form a collection of infinite parallel
%paths perpendicular to the side of $N(P)$, that are \emph{frozen} when
%the slope of the Gibbs measure lies on that boundary of $N(P)$ : with
%probability 1, no dimer will cross these \emph{frozen paths}.
These possibly frozen paths will be the threads of our bead model. When
the slope is not exactly on that boundary of $N(P)$, some dimers may
cross these paths. We will see that these \emph{defects} will play the
role of beads strung along these threads.

Let $\e{e}=(\wv,\bv)$ be an edge of $G$ crossing one of these threads.
When the slope of the Gibbs measure is on the side of
the Newton polygon, this edge appears in the random dimer
configuration with probability 0. In particular, there is no dimer
configuration on the torus $G_1=G/{\ZZ^2}$ containing this edge,
corresponding to a lattice point of
the side of $N(P)$ . As the
cofactor $Q_{\e{e}}(z,w)=Q_{\bv\wv}(z,w)$ is up to a sign the
determinant of the Kasteleyn operator on $G_1\setminus\{\wv,\bv\}$,
it contains only monomials of degree at least 1 in $w$ [otherwise,
there would have been a Gibbs measure with a slope on the side of
$N(P)$ for which $\e{e}$ is a dimer with positive probability, what is
in contradiction with the fact that $\e{e}$ crosses a frozen
path].

We determine now the asymptotic expression for the coupling function
$\Ki_{\gamma,\sm}(\bv_{\xx},\wv)$ corresponding to our magnetic
field $(B_x,B_y)=(c+\beta\gamma\sm,\log\sm)$, between $\wv$
and the black end $\bv_{\xx}$ of a translate
$\e{e}_{\xx}$ of $\e{e}$ by $\xx=(x,y)\in\ZZ^2$
%
%e39 ###
%
\begin{equation}
\Ki_{\gamma,\sm}(\bv_{\xx},\wv)=\iint\frac{z^{-y} w^{x}
Q_{\e{e}}(e^{c+\beta\gamma\sm} z,\sm w)}{P(e^{c+\beta\gamma\sm}
z,\sm w)}\frac{\ud z}{2i\pi z}\frac{\ud w}{2i\pi w}.
\label{eq:Kiperles}
\end{equation}

\begin{prop}\label{Kiperlesasymp}
Denote by $J_{\gamma}(x,\xi)$ the kernel of the bead model.
%
%e40 ###
%
\begin{equation}
J_\gamma(x,\xi)=
\cases{
\displaystyle\int_{[-1,1]}e^{-i\xi
\phi}\bigl(\gamma+i\phi\sqrt{1-\gamma^2}\bigr)^{x}\frac{\ud\phi}{2\pi}, &
\quad if $x\geq0$, \cr
\displaystyle-\int_{\RR\setminus[-1,1]}e^{-i\xi
\phi}\bigl(\gamma+i\phi\sqrt{1-\gamma^2}\bigr)^{x}\frac{\ud\phi}{2\pi},
&\quad if $x<0$.
}
\end{equation}
In the scaling limit
$
\sm\rightarrow0,\quad\quad\sm\beta y\sqrt{1-\gamma^2}\rightarrow
\xi,
$
the coefficients $\Ki_{\gamma,\sm}(\bv_{\xx},\wv)$ have the
following asymptotics:
%
%e41 ###
%
\begin{equation}
\Ki_{\gamma,\sm}(\bv_{\xx},\wv)\sim
\sm\rho_{\e{e}} {\sqrt{1-\gamma^2}}J_\alpha(x,\xi),
\end{equation}
where the quantity $\rho_{\e{e}}$ is given by
%
%e42 ###
%
\begin{equation}
\rho_{\e{e}}=\frac{\partial_2 Q_{\e{e}}(e^c,0)}{\partial_2 P(e^c,0)}.
\end{equation}
$K_{\e{e}}\rho_{\e{e}}$ represents the proportion of this type of
edges among the defects along a thread.
\end{prop}

\begin{pf}
We denote by $f(z,w)$ the rational fraction inside the
integral \eqref{eq:Kiperles}
%
%e43 ###
%
\begin{equation}
f(z,w)=\frac{z^{-y-1} w^{x-1} Q_{\e{e}}(e^{c+\beta\gamma\sm} z,\sm
w)}{P(e^{c+\beta\gamma\sm} z,\sm w)}.
\end{equation}

The integral defining $\Ki_{\gamma,\sm}$ is evaluated by performing
first the integral over $w$.
Suppose first that $x\geq0$. There is no singularity at $w=0$, since
the monomial in $Q$ with lowest degree in $w$ has degree 1. The
only pole in the unit disc is $W_0(z)/\sm$ when $z\in
[e^{-i\theta(\sm)},e^{i\theta(\sm)}]=I_{\sm}$. In this
case, we introduce $\zeta(\phi)=\gamma+i\phi\sqrt{1-\gamma^2}$
with $z=e^{i\beta\sm\phi\sqrt{1-\gamma^2}}$, and get
\begin{eqnarray}
&&\Ki_{\gamma,\sm}(\bv_{\xx},\wv)\nonumber\\
&&\qquad=\int_{I_\sm}\frac{z^{-y}
(W_0(z)/\sm)^{x-1}Q_{\e{e}}(e^{c+\beta\gamma\sm
}z,W_0(z))}{\partial_2
P(z,W_0(z))\sm} \frac{\ud z}{2i\pi z}\nonumber\\[-8pt]\\[-8pt]
&&\qquad=\sqrt{1-\gamma^2}\nonumber\\
&&\quad\qquad{}\times\int_{-1+O(\sm)}^{1+O(\sm)}\tfrac{e^{-i\beta\sm
y\phi\sqrt{1-\gamma^2}}(\zeta(\phi)+O(\sm))^{x-1}
Q_{\e{e}}({ e^{c+\sm\beta\zeta(\phi)}},{\sm\zeta(\phi)+O(\sm
^2)})}{\partial_2
P(e^{c+\sm\beta\zeta(\phi)},\sm\zeta(\phi)+O(\sm^2))}\frac{\beta\ud\phi}{2\pi}.\nonumber
\end{eqnarray}

For a small $\sm$ and a fixed $\phi$, we have
\begin{eqnarray}\label{eq:estimQ}
Q_{\e{e}}\bigl({e^{c+\sm\beta\zeta(\phi)},\sm\zeta(\phi)+O(\sm
^2)}\bigr)&=&\sm\zeta(\phi)\partial_2 Q_{\e{e}} +O(\sm^2),
\\ \label{eq:estimd2P}
\partial_2 P\bigl({ e^{c+\sm\beta\zeta(\phi)},\sm\zeta(\phi)+O(\sm
^2)}\bigr)&=&\partial_2 P+O(\sm),
\end{eqnarray}
where the derivatives of polynomials $P$ and $Q_{\e{e}}$ without
specified point are evaluated at $(e^c,0)$. For the first expansion, we
used the fact that $Q_{\e{e}}(z,0)\equiv0$ and therefore $\partial_1
Q_{\e{e}}=0$.
Therefore,
\begin{eqnarray*}
\Ki_{\gamma,\sm}(\bv_{\xx},\wv)&=&\sm\beta\sqrt{1-\gamma
^2}\frac{\partial_2 Q_{\e{e}}}{\partial_2 P}
\biggl(\int_{[-1,1]}e^{-i\xi\phi}\zeta(\phi)^{x}\frac{\ud\phi
}{2\pi}+O(\sm)\biggr)\\
&=&\sm\beta\sqrt{1-\gamma^2}\rho_e\biggl(\int_{[-1,1]}e^{-i\xi
\phi}(\gamma+i\phi\sqrt{1-\gamma^2})^{x}\frac{\ud\phi}{2\pi
}+O(\sm)\biggr).
\end{eqnarray*}
%
% &=\frac{\sm\beta(e^c\beta\partial_1 Q_{\e{e}}+\partial_2Q_{\e{e}}
%)}{\partial_2 P}
% (\int_{[-\sqrt{1-\gamma^2},\sqrt{1-\gamma}^2]}
% e^{i\xi\phi}(\gamma+i\phi)^{x}\frac{\ud\phi}{2\pi}+O(\sm))

When $x<0$, the rational fraction in the integral has a multiple
pole at $w=0$ which is hard to evaluate directly. However, the
rational fraction is $o(\frac{1}{w})$ as $|w|\rightarrow\infty$,
and hence the sum of all residues in the plane is $0$.

Let us bound the
residues at the simple roots of $P$: $W_1(z),\dots,W_{d-1}(z)$. We
know already from \eqref{eq:P_root_Wj} that there exists $\epsilon$
such that for every
$j\in\{1,\dots,d-1\}$, $|W_j(z)|\geq\epsilon$ for every
$z\in\cercle$. By the same argument of compactness, there
exists a constant $M>0$ such that for $\sm$ small enough,
\begin{eqnarray}
&&\forall j\in\{1,\dots,d-1\}\ \forall
z\in\cercle\nonumber\\[-8pt]\\[-8pt]
&&\qquad|
\partial_2 P(e^{c+\beta\gamma\sm}z,W_j(z)) |\geq
\frac{1}{M}\quad\mbox{and}\quad
\biggl|\frac{Q_{\e{e}}(e^{c+\beta\gamma\sm
}z,W_j(z))}{W_j(z)}\biggr|\leq M\nonumber
\end{eqnarray}
and, therefore,
%
%e44 ###
%
\begin{equation}
\biggl|\frac{(W_j(z)/\sm)^x Q_{\e{e}}(
e^{c+\gamma\beta\sm},W_j(z))}{W_j(z)
\partial_2 P(e^{c+\gamma\beta\sm}z,W_j(z))}\biggr|\leq
(\epsilon/\sm)^{x} M^2=O(\sm^{-x}).
\end{equation}
Thus, the contribution of these residues is negligible as soon as
$x\leq-2$. In that case, we have
\begin{eqnarray}
\Ki_{\gamma,\sm}(\bv_{\xx},\wv)&=&\int_{\cercle}
\Res_{ w=0}f(z,w) \frac{\ud
z}{2i\pi}+\int_{I_\sm}\Res_{ w= W_0(z)/\sm
}f(z,w) \frac{\ud z}{2i\pi}\nonumber\\
&=& -\sum_{j=1}^{d-1}\int_{\cercle} \Res_{
w=W_j(z)/\sm}f(z,w)
\frac{\ud z}{2i\pi}\\
&&{}-\int_{\cercle\setminus I_\sm}\Res_{w= W_0(z)/\sm}f(z,w) \frac{\ud z}{2i\pi}.\nonumber
\end{eqnarray}

Using the estimates \eqref{eq:estimQ} and \eqref{eq:estimd2P}, we
find, using the same change of variable $z=e^{i\beta\sm\sqrt
{1-\gamma^2}\phi}$ as above that
\begin{eqnarray}
&&\int_{\cercle\setminus I_\sm}\Res_{ w= W_0(z)/\sm
}f(z,w)\frac{\ud
z}{2i\pi}\nonumber\\
&&\qquad=\int_{\cercle\setminus I_\sm}\frac{z^{-y}
(W_0(z)/\sm)^{x-1}Q_{\e{e}}(e^{c+\gamma\beta
\sm}z,W_0(z))}{\sm\partial_2
P(e^{c+\gamma\beta\sm}z,W_0(z))}\frac{\ud z}{2i\pi z}\nonumber\\
&&\qquad=\sm\beta\sqrt{1-\gamma^2}\biggl(\int_{
{-\pi/(|\beta|\sm\sqrt{1-\gamma^2})}}^{
-1+o(1)}+\int_{1+o(1)}^{
\pi/(|\beta|\sm\sqrt{1-\gamma^2})}\biggr)
\\
&&\quad\qquad{}\times\frac{e^{-i y\beta\sm\phi\sqrt{1-\gamma^2}} (\zeta(\phi)+O(\sm
))^{x-1} }{\sm\partial_2
P(e^{c+\beta\sm\zeta(\phi)}, \sm\zeta(\phi)+O(\sm^2))}\nonumber\\
&&\quad\qquad{}\times Q_{\e{e}}\bigl(e^{c+\sm\beta\zeta(\phi)}, \sm\zeta(\phi)+O(\sm
^2)\bigr)\frac{\ud\phi}{2\pi}.\nonumber
\end{eqnarray}
%
%_{-\pi/|\beta|}^{-\sqrt{1-\gamma^2}+O(\sm)}+\int_{\sqrt{1-\gamma^2}+O(

By Lebesgue's dominated convergence theorem and by \eqref{eq:estimQ}
and \eqref{eq:estimd2P}, the integral is asymptotic in the scaling limit
$\sm\rightarrow0$, $\sm\beta y\sqrt{1-\gamma^2}\rightarrow\xi$
to the following expression:
%
%e45 ###
%
\begin{equation}
\frac{\partial_2Q_{\e{e}}}{\partial_2 P}
\int_{\RR\setminus[-1,1]}       e^{-i\phi
\xi} \bigl(\gamma+i\phi\sqrt{1-\gamma^2}\bigr)^x \frac{\ud\phi}{2\pi}
\end{equation}
and thus,
%
%e46 ###
%
\begin{eqnarray}
\Ki_{\gamma,\sm}(\bv_{\xx},\wv)&=&\frac{\sm\beta\sqrt{1-\gamma
^2}\partial_2 Q_{\e{e}}}{\partial_2 P}\nonumber\\[-8pt]\\[-8pt]
&&{}\times\biggl(-\int_{\RR\setminus[-1,1]}
e^{-i\phi\xi} \bigl(\gamma+i\phi\sqrt{1-\gamma^2}\bigr)^x \frac{\ud\phi
}{2\pi}+o(1)\biggr).\nonumber
\end{eqnarray}

When $x=-1$, the residues at the poles $W_1(z),\dots,W_{d-1}(z)$ are
not negligible any more. However, in this case, the pole at $w=0$ is
simple. A direct evaluation of the integral shows
\begin{eqnarray}
\Ki_{\gamma,\sm}(\bv_{(-1,y)},\wv)&=&\int_{\cercle}\Res
_{ w=0}f(z,w)\frac{\ud z}{2i\pi} + \int_{I_\sm
}\Res_{ w=W_0(z)/\sm} f(z,w)\frac{\ud
z}{2i\pi}\nonumber\\
&=&\int_{\cercle}\frac{z^{-y}\sm\partial_2 Q_{\e{e}}(e^{c+\beta
\gamma\sm}z,0)}{P(e^{c+\beta\gamma\sm}z,0)}\frac{\ud z}{2i\pi z}\\
&&{}+\int_{I_\sm}\frac{z^{-y}(W_0(z)/\sm
)^{-2}Q_{\e{e}}(e^{c+\beta\gamma\sm}z,W_0(z))}{\sm\partial_2
P(e^{c+\beta\gamma\sm} z,W_0(z))}\frac{\ud z}{2i\pi z}.\nonumber
\end{eqnarray}

Posing again $z=e^{i\beta\sm\phi\sqrt{1-\gamma^2}}$ and $\zeta
(\phi)=\gamma+i\phi\sqrt{1-\gamma^2}$, one has
\begin{eqnarray}
\Ki_{\gamma,\sm}(\bv_{\xx},\wv)&=& \sm^2\beta\sqrt{1-\gamma
^2}\nonumber\\
&&{}\times\int_{-\pi/(|\beta|\sm\sqrt{1-\gamma^2})}^{\pi
/(|\beta|\sm\sqrt{1-\gamma^2})} \frac{e^{-i\beta\varphi\sm
y\sqrt{1-\gamma^2}} \partial_2 Q_{\e{e}}(e^{c+\beta\sm\zeta(\phi
)},0)} {P(e^{c+\beta\sm\zeta(\phi)},0)} \frac{\ud\phi} {2\pi}\\
 %\nonumber\\[-8pt]\\[-8pt]
&&{}+\sm^2\beta\sqrt{1-\gamma^2}\int_{
-1+o(1)}^{1+o(1)} \frac{e^{-i\beta\varphi\sm
y\sqrt{1-\gamma^2}}Q_{\e{e}}(e^{c+\beta\sm\zeta(\phi
)},0)}{W_0(e^{i\beta\phi\sm\sqrt{1-\gamma^2}})^2\partial_2
P(e^{c+\beta\sm\zeta(\phi)},0)}\frac{\ud\phi}{2\pi}.\nonumber
\end{eqnarray}
%
% what can be rewritten as
% \begin{multline}
% \Ki_{\gamma,\sm}(\bv_{\xx},\wv)= \sm^2\beta(\int_{[-\pi/|\beta|,
% \end{multline}

Using the following estimates from Taylor's formula
%
%e47 ###
%
\begin{eqnarray}
\qquad \partial_2 Q\bigl(e^{c+\beta\sm\zeta(\phi)},0\bigr) &=& \partial_2 Q_{\e
{e}}+O(\sm),\\
P\bigl(e^{c+\beta\sm\zeta(\phi)},0\bigr)&=&P\bigl(e^{c+\beta\gamma\sm+i\beta\sm
\zeta(\phi)},W_0(z)\bigr)\nonumber\\
&&{}-W_0(z)\partial_2 P\bigl(e^{c+\beta\sm\zeta(\phi
)},W_0(z)\bigr)+O(\sm^2)\\
&=&0-\sm\zeta(\phi)\bigl(\partial_2 P+O(\sm)\bigr)=-\sm\zeta(\phi)\partial
_2 P +O(\sm^2)\nonumber
\end{eqnarray}
together with \eqref{eq:estimQ} and \eqref{eq:estimd2P}, and applying
Lebesgue dominated convergence theorem after an integration by parts,
one can prove that in the scaling limit, the coefficient of the inverse
Kasteleyn operator $\Ki_{\gamma,\sm}(\bv_{\xx},\wv)$ is
asymptotic to
\[
-\sm\beta\sqrt{1-\gamma^2}\frac{\partial_2 Q_{\e{e}}}{\partial
_2 P} \int_{\RR\setminus[-1,1]}e^{-i\phi\xi}\bigl(\gamma+i\phi\sqrt
{1-\gamma^2}\bigr)^{-1}\frac{\ud\phi}{2\pi}.
\]
\upqed
\end{pf}

\begin{rem}
The ratio
$
{ \partial_2 Q_{\e{e}}}/ {\partial_2 P}
$
controls the density of the copies of $\e{e}$ in the limiting bead
model. If one plugs the value $\phi=1$ into \eqref{eq:estimQ} and
\eqref{eq:estimd2P}, one can see that it is up to terms of higher
order in $\sm$, equal to
%
%e48 ###
%
\begin{equation}
\frac{i Q_{\e{e}}(e^{c+\beta\gamma\sm}z_0,\sm w_0)}{i\partial_2
P(e^{c+\beta\gamma\sm}z_0,\sm w_0)\sm w_0},
\end{equation}
where $(z_0,w_0)$ are zeros of the characteristic polynomial
$ P(e^{c+\beta\gamma\sm}\cdot,\sm\cdot)$
on the unit torus. When multiplied by $\K_{\e{e}}$, the numerator is
the length
of the dual edge $\e{e}^*$ in the natural application from
the dual graph $G^*$ to $\RR^2$ described below while the denominator
is that of the vertical side of the fundamental domain.

\begin{lem}[(\cite{Boutdensities})]\label{lem:mapping}
Let $(e^{B_x}z_0,e^{B_y}w_0)$ be a root of the characteristic
polynomial, with $(z_0,w_0)$ on the unit torus. The 1-form
%
%e49 ###
%
\begin{equation}
\e{e}=(\wv,\bv)\mapsto i K_{\wv\bv}(e^{B_x}z_0,e^{B_y}w_0) Q_{\bv
\wv}(e^{B_x}z_0,e^{B_y}w_0)
\end{equation}
is a divergence-free flow. Its dual is therefore the gradient of a
mapping from $G^*$ to $\RR^2\simeq\CC$.

This mapping $\Psi^*$ is $\ZZ^2$-periodic and the symmetries of its
range are generated by
\[
\hat{\mathrm{x}}=i e^{B_x}z_0\partial_1
P(e^{B_x}z_0,e^{B_y}w_0)\quad \mbox{and}\quad \hat{\mathrm{y}}=i
e^{B_y}w_0\partial_2 P(e^{B_x}z_0,e^{B_y}w_0).
\]
%
% There exists a map $\Psi$ from $G^*$ to $\mathbbm{R}^2\simeq
% dual of the edge $\e{e}=(\wv,\bv)$ is represented by the complex
%number
% $\K(\wv,\bv) Q(z_0,w_0)$.
\end{lem}

This application $\Psi$ coincide with the notion of isoradial
embedding for dimer models with critical weights (see \cite{KeCrit})
and gives a geometry well adapted for the study of liquid measures on
dimer configurations (see, e.g.,~\cite{Boutdensities}).
\end{rem}

Proposition \ref{Kiperlesasymp} shows that the kernel giving the
correlations has the same form as the original bead model.
However, in order to recover fully the bead model, one can not just
look at one type of edges on threads but at all of them.

Since the frozen paths have been chosen to cross no dimer when the
slope is on the side of the particular Newton
polygon we are looking at, they are bordered by white vertices on their
left and black vertices on their right.
For a reason of parity between white and black vertices, there is no
dimer configuration of the graph $G_1$ deprived of the projection of
these two vertices having a height change\footnote{The term \emph
{height change} here is an abuse of notations, since the difference
between the reference unit flow $\conf_0$ and the one corresponding to
any dimer configuration of $G_1$ deprived of the two vertices has a
nonzero divergence. What we mean here by this expression is, in fact,
the powers in $z$ and $w$ in the weight of the configuration computed
using the magnetically altered Kasteleyn operator $K(z,w)$ divided by
that of the reference dimer configuration.} on the side of $N(P)$.
Therefore, the arguments of the proof of Proposition \ref
{Kiperlesasymp} can be applied to obtain similar asymptotics as those
given in that proposition for the coefficients of $\Ki_{\gamma,\sm}$
between these vertices.
\begin{prop} Let $\bv_{x,y}$ and $\wv$ be respectively a black and a
white vertex each bordering one of these paths, and in
fundamental domains separated by a lattice translation $(x,y)$.
In the scaling limit,
%
%e50 ###
%
\begin{equation}
\sm\rightarrow0,\qquad\sm\beta y\sqrt{1-\gamma^2} \rightarrow\xi,
\end{equation}
the coefficient $\Ki_{\gamma,\sm}(\bv_{x,y},\wv)$ has the
following asymptotics
%
%e51 ###
%
\begin{equation}
\Ki_{\gamma,\sm}(\bv_{x,y},\wv) \sim\sm\rho_{\bv\wv}J_{\gamma
}(x,\xi),
\label{eq:asympperles2}
\end{equation}
where
%
%e52 ###
%
\begin{equation}
\rho_{\bv\wv}=\frac{\partial_2 Q_{\bv\wv}(e^c,0)}{\partial_2 P(e^c,0)}.
\end{equation}
\end{prop}

These coefficients $\rho_{\bv\wv}$ are in fact the product of two
terms, one depending only on $\bv$ and the other on $\wv$. This
property is stated in the following lemma.
\begin{lem}\label{lem:Qrank1}
The rank of the matrix $\partial_2 Q(e^c,0)$, restricted to
projections of vertices bordering a thread is equal to 1. In
particular, for any $\bv$ and any $\wv$ bordering a thread, there
exist $U_{\bv}$ and $V_{\wv}$ such that
%
%e53 ###
%
\begin{equation}
\rho_{\bv\wv}=U_{\bv}V_{\wv}.
\end{equation}
\end{lem}

\begin{pf}
The matrix $\partial_2 Q(e^c,0)$ is the limit when $\sm$ goes to zero of
%
%e54 ###
%
\begin{equation}
\frac{1}{\sm w_0}Q(e^{c+\beta\gamma\sm}z_0,\sm w_0),
\end{equation}
where $(z_0,w_0)$ are zeros on the unit torus of $P(e^{c+\beta\gamma
\sm}z,\sm w)$. These zeros depend on $\sm$ and $\gamma$ and their
the first term in their expansion in $\sm$ is obtained by plugging
$\phi=\pm1$ into equations \eqref{eq:estimQ} and \eqref{eq:estimd2P}:
%
%e55 ###
%
\begin{equation}
z_0=1+O(\sm),\qquad w_0= \gamma\pm i\sqrt{1-\gamma^2}+O(\sm).
\end{equation}
Since $(e^{c+\beta\gamma\sm}z_0,\sm w_0)$ is not real, it is a
simple zero of $P=\det K$, since the mapping from the spectral curve to
its am{\oe}ba is 2-to-1 out of the real locus. $Q$ is the comatrix of
$K$, its rank is 1 at a simple zero of $P$. Therefore, as the limit of
sequence of rank 1 matrices, the matrix
$\partial_2 Q(e^c,0)$ has a rank a most 1. As there is at least a
nonzero entry in this matrix, the rank is equal to 1.

The coefficient $\rho_{\bv\wv}$ is a multiple of the entry $(\bv
,\wv)$ of this matrix.
Its decomposition into a product comes from the representation of a
rank-1 matrix as a tensor product of a vector and a linear form.
\end{pf}

%
% As before, the coefficient is the ratio between some derivative of
%the entry $Q_{\bv\wv}$ and one of the polynomial $P$. More precisely,
%
% Since $Q(z,w)$ is the comatrix of $K(z,w)$, its rank at a simple zero
%of $P(z,w)=\det K(z,w)$ is equal to 1. Therefore,theis $(z_0,w_0)$ is
%a zero of $P(e^{c+\beta\gamma\sm}z,\sm w)=\det K(e^{c+\beta\gamma\sm}z,

%s5.4 ###
\subsection{Convergence to the bead model}

We already said that the threads of our bead model would be the
infinite collection of vertical paths, translated one from another,
that are frozen, \textit{that is,} they do not cross any dimer when the
slope of the measure lies on the boundary of the Newton polygon. The
beads are represented by the dimers crossing these paths when the
magnetic field lies in one of the vertical tentacles of the am{\oe}ba.

Like in Section~\ref{sec:bead_hex}, as the magnetic field goes deeper
into the tentacle of the am{\oe}ba, the picture of the graph in the
plane is rescaled in such a way that although the probability of seeing
a particular dimer crossing these ``almost-frozen'' paths goes to zero,
the average number of such edges by centimeter of thread stays almost
constant. The scaling limit we perform is
%
%e56 ###
%
\begin{equation}
\sm\rightarrow0,\qquad\sm\beta y\sqrt{1-\gamma^2}\rightarrow\xi
\in\RR.
\end{equation}

To find the limiting distribution of this beads, we first evaluate the
quantities
%
%e57 ###
%
\begin{equation}
\mathbb{E}\biggl[\frac{X_{I_1} !}{(X_{I_1}-n_1)!}\cdots\frac{X_{I_k}
!}{(X_{I_k}-n_k)!}\biggr],
\end{equation}
where $X_{I_j}$ is the number of dimers crossing the (rescaled) thread
interval. We look in detail at the case $k=1$ when only one thread
interval is at stake. The other cases are similar. For a given $n$, and
a fixed value of $\gamma$ and of the scaling parameter $\sm$, we have
%
%e58 ###
%
\begin{equation}
\mathbb{E}_{\gamma,\sm}\biggl[\frac{X_I!}{(X_I-n)!}\biggr]=\mathop{\sum_{\e
{e}_1,\dots,\e{e}_n\in I }}_{ \mathrm{distinct}} \probind{\gamma,\sm
}{\e{e}_1,\dots,\e{e}_n\in\conf},
\end{equation}
where the sum is over all possible $n$-tuples of edges crossing the
thread interval $I$. The edges crossing $I$ are labeled by their type
(i.e., their projection on $G_1=G/\ZZ^2$) and the coordinates of
the fundamental domain they belong to. The edge $\e{e}_{\xx}^j$
represents the edge of type $j$ in the fundamental domain with
coordinates $\xx=(x,y)$. The type label $j$ ranges from 1 to $d$.
Since the probability of having two such edges in the same fundamental
domain is negligible, we can rewrite this sum of probabilities, as a
sum over the fundamental domains and the types of edges crossing the
thread interval.
%
%e59 ###
%
\begin{equation}\label{eq:perlesgen}
\mathbb{E}_{\gamma,\sm}\biggl[\frac{X_I!}{(X_I-n)!}\biggr]=\mathop{\sum_{
{\xx
_1,\dots,\xx_n\in I }}}_{ \mathrm{distinct}} \sum_{j_1,\dots,j_n=1}^n
\probind{\gamma,\sm}{\e{e}_{\xx_1}^{j_1},\dots,\e{e}_{\xx
_n}^{j_n}\in\conf} +O(\sm).
\end{equation}
The different probabilities $\probind{\gamma,\sm}{\e{e}_{\xx
_1}^{j_1},\dots,\e{e}_{\xx_n}^{j_n}\in\conf}$ are given by the
determinant of a $n\times n$ matrix that is equal to
%
%e60 ###
%
\begin{eqnarray}\label{eq:detperles}
&&\bigl(\sm\beta\sqrt{1-\gamma^2}\bigr)^n \Biggl(\prod_{l=1}^n \K_{j_l}
\Biggr)\nonumber\\[-8pt]\\[-8pt]
&&\qquad{} \times\det_{1\leq k,l\leq n} \bigl[\rho_{j_k j_l}J_\gamma\bigl(x_k-x_l,\sm
\beta\sqrt{1-\gamma^2}(y_k-y_l)\bigr)\bigr]+O(\sm^{n+1}).\nonumber
\end{eqnarray}
Since $\rho_{jk}$ is the product of two terms $U_j V_k$, we can carry
them out of the determinant by $n$-linearity these coefficients,
equation \eqref{eq:detperles} becomes
\begin{eqnarray}
&&\bigl(\sm\beta\sqrt{1-\gamma^2}\bigr)^n \Biggl(\prod_{l=1}^n \K_{j_l}
\Biggr)\nonumber\\
&&\quad{} \times\det_{1\leq k,l\leq n} \bigl[U_{j_k}V_{j_l}J_\gamma
\bigl(x_k-x_l,\sm\beta\sqrt{1-\gamma^2}(y_k-y_l)\bigr)\bigr]+O(\sm^{n+1})
\nonumber\\
&&\qquad=\bigl(\sm\beta\sqrt{1-\gamma^2}\bigr)^n \Biggl(\prod_{l=1}^n \K
_{j_l}U_{j_l}V_{j_l}\Biggr)\nonumber\\[-8pt]\\[-8pt]
&&\quad\qquad{} \times\det_{1\leq k,l\leq n}
\bigl[J_\gamma\bigl(x_k-x_l,\sm\beta\sqrt{1-\gamma^2}(y_k-y_l)\bigr)\bigr]+O(\sm
^{n+1})\nonumber \\
&&\qquad=\bigl(\sm\beta\sqrt{1-\gamma^2}\bigr)^n \Biggl(\prod_{l=1}^n \K_{j_l}\rho
_{j_l}\Biggr)\nonumber\\
&&\quad\qquad{} \times\det_{1\leq k,l\leq n} \bigl[J_\gamma
\bigl(x_k-x_l,\sm\beta\sqrt{1-\gamma^2}(y_k-y_l)\bigr)\bigr]+O(\sm^{n+1}).\nonumber
\end{eqnarray}
Summing now over the different types of edges, one gets
\begin{eqnarray}
&&\sum_{j_1,\dots,j_n=1}^{d} \probind{\gamma,\sm}{\e{e}_{\xx
_1}^{j_1},\dots,\e{e}_{\xx_n}^{j_n}\in\conf} \nonumber\\
&&\qquad=\bigl(\sm\beta\sqrt{1-\gamma^2}\bigr)^n \Biggl(\sum_{j=1}^{d} \K_{j_l}\rho
_{j_l}\Biggr)^n \\
&&\quad\qquad{} \times\det_{1\leq k,l\leq n} \bigl[J_\gamma\bigl(x_k-x_l,\sm
\beta\sqrt{1-\gamma^2}(y_k-y_l)\bigr)\bigr]+O(\sm^{n+1}).\nonumber
\end{eqnarray}
$\K_j \rho_j$ is the proportion of edges of type $j$ along the
thread. These coefficients sum up to 1, and we have finally that
expression \eqref{eq:perlesgen} is a Riemann sum for
%
%e61 ###
%
\begin{equation}
\int\cdots\int_{I^n} \det_{1\leq k,l\leq n}[J_{\gamma}(x_k-x_l,\xi
_k-\xi_l)]\ud^n\xi
\end{equation}
and the same argument of domination as in the proof of Theorem \ref
{thm:bead1} implies that the distribution of $X_I$ converges to the
distribution of beads in the interval $I$ in the bead model of
parameter $\gamma$.
The generalization to any finite dimensional distribution is
notationally cumbersome, but straightforward. These considerations give
thus the proof of the following theorem.
\begin{thm} \label{thm:bead2}
Let $\gamma\in(-1,1)$. In the scaling limit $\sm\rightarrow0, \sm
\beta y\sqrt{1-\gamma^2}\rightarrow\xi$,
the point process describing the position of rare edges on the threads,
identified with the almost frozen paths converges to the \emph{bead
model} of index $\gamma$, \textit{that is,} the determinantal point
process on $\ZZ\times\RR$ with kernel $J_{\gamma}$.
% In particular, the marginal along each thread is the GUE process.
\end{thm}

Recall that $\gamma$ describes the different possible ways to go deep
into a tentacle. This theorem states that the bead model, with its
1-parameter family of Gibbs measures, is the universal limiting
behavior of any dimer model on a bipartite periodic planar graph when
the order parameters $(B_x,B_y)$ go to infinity in staying in the
liquid phase.

%s6 ###
\section{Interaction between bead models}

It often happens that a side of the Newton polygon is not the result of
a unique frozen path, but that different paths give the same constraint
on the slope. In that case, we do not have just one family of frozen
paths, but several parallel families of thread, carrying all in the
scaling limit a bead model. In this section, we describe the
interaction between these different bead models in the case of the
generic case of the honeycomb lattice $H$ with a $n\times m$
fundamental domain.

\begin{figure}

\includegraphics{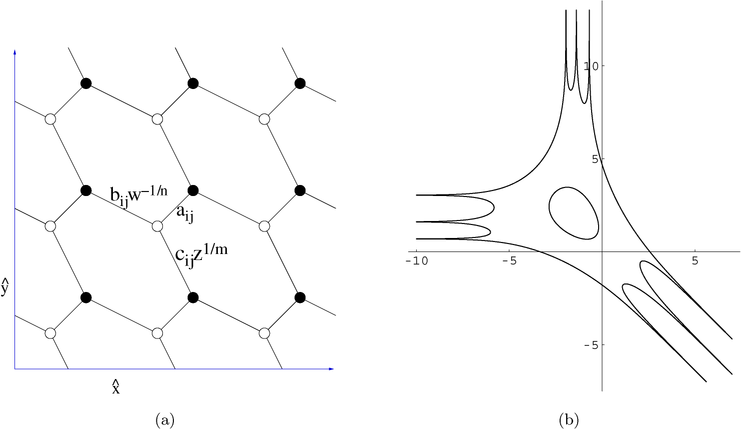}

\caption{\textup{(a)} A $3\times3$ fundamental domain of the
honeycomb lattice. Weights of edges around white vertex $\wv_{ij}$ are
$a_{ij}$, $b_{ij}$ and $c_{ij}$. The Fourier multipliers have been
distributed over the edges by gauge transformation. \textup{(b)} The am{\oe
}ba for generic weights on the graph represented on the left panel.
This am{\oe}ba presents a gaseous phase, and three tentacles in three
directions, each corresponding to a collection of frozen paths drawn on
Figure~\textup{\protect\ref{fig:hexgenericpaths}}.}
\label{fig:bare_hexgeneric}
\end{figure}

The fundamental domain of this periodic planar graph is represented on
Figure \ref{fig:bare_hexgeneric} for $n=m=3$.
The vertices of the fundamental domain are labeled by two integers, $i$
and $j$ ranging from 1 to $n$, and from $1$ to $m,$ respectively. The
weights of the edges around the white vertex labeled by $(i,j)$ are
denoted by $a_{ij}$, $b_{ij}$ and $c_{ij}$. By an appropriate gauge
transformation~\cite{KOS}, one can spread the factors $z$ and $w$ in
the magnetically altered Kasteleyn matrix $K(z,w)$ so that the
coefficients of this operator are $a_{ij}$, $b_{ij}w^{-1/n}$ and
$c_{ij}z^{1/m}$.
The reference dimer configuration we will use is the configuration
containing all the $a$-edges.

%f8 ###
%

%f9 ###
%
\begin{figure}[b]

\includegraphics{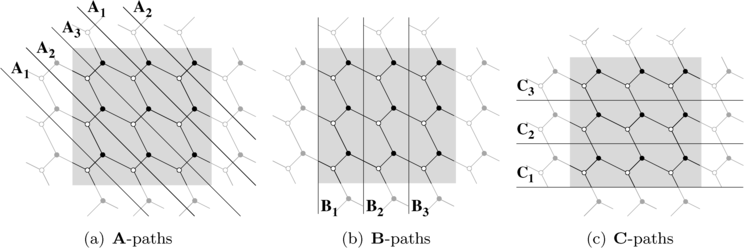}

\caption{The three classes of possible frozen paths:
$\mathbf{A}$, $\mathbf{B}$ and $\mathbf{C}$.}
\label{fig:hexgenericpaths}
\end{figure}

One distinguishes three special classes of dual cycles in $G_1$, say
$\mathbf{A}$, $\mathbf{B}$ and $\mathbf{C}$ that cross only edges of
a given type (resp. $a$, $b$ and $c$). The lifts of these
classes to $H$ are represented on Figure~\ref{fig:hexgenericpaths}.
The $\mathbf{B}$ class is constituted by $n$ vertical straight paths
with homology class $(0,1)$, whereas the $\mathbf{C}$ class contains
the $m$ horizontal straight paths with homology class $(1,0)$. The
$\mathbf{A}$ class contains $d=\gcd(n,m)$ paths with homology $(\frac
{m}{d},\frac{n}{d})$.
The three classes of cycles lift to $H$, forming three classes of
parallel families of straight lines.

The Newton polygon of the weighted dimer graph $G$ is a right triangle.
Each side of the triangle corresponds to Gibbs measures for which all
the paths of one of the three classes are frozen.
The horizontal side contains $n+1$ lattice points. The am{\oe}ba of
the associated spectral curve exhibits $n$ vertical tentacles
separating $n+1$ unbounded complementary components---solid phases in
the phase diagram---each corresponding to a lattice point of the
horizontal side.

The frozen configurations are obtained as follows: When the Gibbs
measure's slope lies, say, on the horizontal side of the Newton
polygon, there is almost surely no edge crossing the $\mathbf{B}$
class of paths. These paths delimit thin strips where one sees either
an infinite succession of $a$-edges or an infinite collection of
$c$-edges. On $G_1$ and in presence of a magnetic field $(B_x,B_y)$,
the associated patterns have weights $\prod_{i=1}^m a_{ij}$ and
$e^{B_x}\prod_{i=1}^m c_{ij}$. The patterns with the highest weight
correspond to the type of columns appearing in the configuration with
probability 1. When the horizontal component of the magnetic field is
very negative, the weight of the ``$a$'' patterns are greater than the
``$c$'' ones, but as $B_x$ increases, the weight of the second pattern
becomes more important, and at some point, it becomes bigger than the
first one. In the graph $G$, the $a$-edges that were filling the space
between the two $\mathbf{B}$ paths switch to $c$-edges.
Generically, the values of $B_x$ corresponding to such a switch are all
different. They correspond to the absciss{\ae} of the vertical
tentacles of the amo{\oe}ba.

In a fixed window and when $B_y$ is very large, one sees columns of
edges of the same type ($a$, or $c$) with a probability close to 1.
When the magnetic field lies in a tentacle of the phase diagram, the
system hesitates between two states for a given type of columns. With a
probability $p$ bounded away from 0 and 1, the column is filled with
$a$-edges, and with probability $1-p$, it is filled with $c$-edges.
If one rescales vertically the graph in the same time as $B_y$ goes to
$+\infty$, then one will be able to see the transition between these
two possibilities: Between the two types of fillings, a $b$-edge is
inserted. The edge creates a defect in the neighboring column that is
supposed to be frozen. This discussion is quantified in the following
proposition.

\begin{prop}\label{prop:edge01}
If
%
%e62 ###
%
\begin{equation}
B_x<\sum_{i=1}^m\log\frac{a_{ij_0} } {c_{ij_0}},
\end{equation}
then with probability going to 1 when $B_y$ goes to infinity, the
columns of type $j_0$ will be filled with $a$-edges.

If the inequality is strict in the other direction,
then they will be filled with $c$-edges with probability going to 1
when $B_y$ goes to infinity.
\end{prop}

\begin{pf}
When the vertical component of the magnetic field is very negative, the
main contribution to the characteristic polynomial is given by the
configurations on $G_1$ that contain no $b$-edges. One can choose to
fill each strip of $G_1$ between two consecutive frozen cycles either
by $a$-edges or by $c$-edges. Choosing a filling with $c$ edges induces
a height change of $m$ in the vertical direction. Therefore,
%
%e63 ###
%
\begin{equation}
\qquad P^B(z,w)=P(e^{B_x}z,e^{B_y}w)=\prod_{j=1}^n\Biggl(\prod_{i=1}^m
a_{ij}-(-1)^m e^{B_x}z\prod_{i=1}^m c_{ij}\Biggr)+O(\sm),
\end{equation}
where $t=e^{B_y}$ is small.

Let $\bv$ and $\wv$ two vertices in the same strip $j_0$ in $G_1$.
Denote by $\bv_{y}$ and $\wv_0$ lifts of $\bv$ and $\wv$ in the
same column of $G$, separated by $y$ fundamental domains.
%If $\wv$ and $\bv$ are two vertices in the same strip, t
The entry of the inverse Kasteleyn operator between these two vertices
is easily evaluated: recall that $Q_{\bv\wv}$ is the characteristic
polynomial of the graph $G$ where all the translated of $\bv$ and $\wv
$, as well as all the edges connected to these vertices, have been
removed. Repeating the argument given above, we find the main
contribution to it,
%
%e64 ###
%
\begin{eqnarray}
Q^B_{\bv\wv}(z,w)&=&Q_{\bv\wv}(e^{B_x}z,e^{B_y}w)\nonumber\\[-8pt]\\[-8pt]
&=&M_{\bv\wv
}z^\delta\mathop{\prod_{j=1}}_{ j\neq j_0}^n\Biggl(\prod_{i=1}^m
a_{ij}-(-1)^m e^{B_x}z\prod_{i=1}^m c_{ij}\Biggr)+O(\sm),\nonumber
\end{eqnarray}
where $M_{\bv\wv}z^\delta$ is the weight of the dimer configuration
of the strip $j_0$ of $G_1$ deprived of $\bv$ and $\wv$.
The coefficient of the inverse Kasteleyn operator corresponding to
these two vertices whose fundamental domains are separated by the
lattice vector $(x,y)$ is given by
\begin{eqnarray}\label{eq:Ki01}
\qquad\Ki_B(\bv_y,\wv_0)&=&\iint_{\tore} \frac{z^{-y}w^{0} Q^B_{\bv\wv
}(z,w)}{P^B(z,w)}\frac{\ud z}{2i\pi z} \frac{\ud w}{2i\pi w} \\
&=&\iint_{\tore} \frac{z^{-y+\delta} M_{\bv\wv}}{(\prod
_{i=1}^m a_{ij_0}-(-1)^m e^{B_x}z\prod_{i=1}^m c_{ij_0})}\frac
{\ud z}{2i\pi z} \frac{\ud w}{2i\pi w}\nonumber\\[-8pt]\\[-8pt]
&&{}+O(\sm) \nonumber\\
&=&\int_{\cercle}\frac{z^{-y+\delta} M_{\bv\wv}}{(\prod
_{i=1}^m a_{ij_0}-(-1)^m e^{B_x}z\prod_{i=1}^m c_{ij_0})}\frac
{\ud z}{2i\pi z}+O(\sm).
\end{eqnarray}

Suppose that $B_x<\sum_{i=1}^m\log\frac{a_{ij_0} } {c_{ij_0}}$, the
other case is similar.
In that case, the pole located at
%
%e65 ###
%
\begin{equation}
z=(-1)^m \frac{\prod_{i=1}^m a_{ij_0}}{e^{B_x}\prod_{i=1}^m c_{ij_0}}
\end{equation}
is outside of the unit disk. If $y-\delta<0$, then there is no pole at
all in the unit disk and, therefore, the integral over $z$ is zero.
However, if $y-\delta\geq0$, then the integral equals the opposite of
the residue at $z_0$, and
%
%e66 ###
%
\begin{equation}
\qquad \Ki_B(\bv_y,\wv_0)=\frac{M_{\bv\wv}}{(-1)^m e^{B_x}\prod_{i=1}^m
c_{ij_0}}\biggl(\frac{\prod_{i=1}^m a_{ij_0}}{(-1)^m e^{B_x}\prod
_{i=1}^m c_{ij_0}}\biggr)^{-y-1}+O(\sm).
\end{equation}
When $\bv_y$ and $\wv_0$ are the ends of an edge with weight
$a_{i_0j_0}$, $M_{\bv\wv}=\prod_{i\neq i_0} a_{ij_0}$ and $y$ and
$\delta$ equal 0. It follows that
%
%e67 ###
%
\begin{equation}
\Ki_B(\bv_y,\wv_0)=\frac{M_{\bv\wv}}{\prod_{i=1}^m
a_{ij_0}}+O(\sm)=\frac{1}{a_{i_0j_0}}+O(\sm).
\end{equation}
Thus, the probability of this $a$-edge, given by $a_{i_0j_0} \Ki_B(\bv
_y,\wv_0)$ goes to 1 when $\sm$ goes to zero.

On the other hand, in the case when $\bv$ and $\wv$ are the ends of a
``$c$''-type edge, then either $y=-1$ or $\delta=1$. In both cases, $\Ki
_B(\bv_y,\wv_0)$ is $O(\sm)$. Thus, the probability of this edge
goes to zero when $\sm$ goes to zero. Such an edge is called \emph
{nontypical}.%\rightqed
\end{pf}

A sequence of nontypical edges in a frozen column is initiated by the
presence of a bead (a ``$b$''-edge) crossing a neighboring wire.
The analysis we made of the inverse Kasteleyn operator allows us to
determine the distribution of the length of the sequence of nontypical edges in a
frozen column.

\begin{prop}\label{prop:nontypicaledges}
The length of a succession of nontypical edges in a frozen column has
a geometric distribution in the limit. The parameter of the geometric
distribution has an explicit expression in terms of ratios of lengths
of dual edges.
\end{prop}

\begin{pf}
We suppose that in the frozen column $j_0$, we only see $a$-edges with
probability close to 1.
The  inequality
%
%e68 ###
%
\begin{equation}
B_x<\sum_{i=1}^m \log\frac{a_{ij_0}}{c_{ij_0}}
\end{equation}
is satisfied. Since we work only in one column, we will drop the index
$j_0$ for the sake of simplicity. See Figure~\ref{fig:illustprop} for
an illustration of the notations. Given that the edge $\e{e}=(\wv,\bv
)$ with weight $b_{i_0}$ is present in the dimer configuration, we
compute the probability of seeing $N$ successive $c$-edges after this
bead. Denote by $\e{e}^c_1=(\wv_1,\bv_1),\dots,\e{e}^c_N(\wv
_N,\bv_N)$ the $N$ $c$-edges. The weights of the edges around vertex
$\wv_i$ are $a_{[i]},b_{[i]}, c_{[i]}$, where $[i]=(i_0+i \operatorname{mod}
m)+1$. The conditional probability we want to compute is given by
the following formula:
\begin{eqnarray}
&&\prob{\e{e}\in\conf\mbox{ and }\forall i=1,\dots,N \e{e}^c_i\in
\conf| \e{e}\in\conf}\nonumber\\
&&\qquad=\frac{\prob{\e{e}\in\conf\mbox{ and }\forall i=1,\dots,N \e
{e}^c_i\in\conf}} {\prob{\e{e}\in\conf}}\\
&&\qquad=\Biggl(\prod_{i=1}^N
c_{[i]}\Biggr)\frac{\det A_{N+1} }{\Ki(\bv,\wv)},\nonumber
\end{eqnarray}
where $A_{N+1}$ is the following square matrix whose entries are
inverse Kasteleyn operator coefficients:
%
%e69 ###
%
\begin{equation}
A=\left[
\begin{array}{c|ccc}
\Ki(\bv,\wv) &\cdots&\Ki(\bv_j,\wv) &\cdots\\
\hline\vdots&\ddots&\vdots& \\
\Ki(\bv,\wv_i) & &\Ki(\bv_j,\wv_i) & \\
\vdots& &\vdots&\ddots
\end{array}
\right].
\end{equation}
Since $\wv$ and the white vertices $\wv_i$ stand on different sides
of a frozen path, the associated coefficients $\Ki(\bv,\wv_i)$ are
$O(\sm)$. More precisely, from \eqref{eq:asympperles2}, one has
%
%e70 ###
%
\begin{equation}
\Ki(\bv,\wv_i)=\sm\rho_{\bv\wv_i}+O(t^2).
\end{equation}
Besides, if $i\leq j$, the power of $z$ in the numerator of \eqref
{eq:Ki01} is positive, and it follows from computations made above that
$\Ki(\bv_j,\wv_i)$ is also $O(\sm)$.

%f10 ###
%
\begin{figure}

\includegraphics{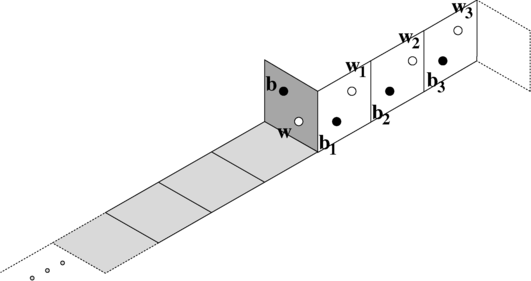}

\caption{Illustration of Proposition \textup{\protect\ref
{prop:nontypicaledges}}. Here is represented a frozen column of
(horizontal) $a$-rhombi (in light grey), perturbed by the presence of a
defect (the $b$-rhombus corresponding to the edge $(\wv,\bv)$),
followed by a finite sequence of $c$-rhombi.}\label{fig:illustprop}
\end{figure}

For $i\in\{1,\dots,N\}$ and with the convention that $\wv_0=\wv$,
$\wv_{i-1}$ and $\bv_i$ are the ends of an $a$-edge with weight
$a_{[i-1]}$. The same computations as above show that $\Ki(\bv_i,\wv
_{i-1})=\frac{1}{a_{[i-1]}}$.
As a consequence, the asymptotic expansion of the determinant of
$A_{N+1}$ is given by the product of these elements just above the
diagonal times the last element of the first column.
%
%e71 ###
%
\begin{equation}
\det A_{N+1}= \sm\rho_{\bv\wv_N}\prod_{i=1}^{N} \frac
{1}{a_{[i-1]}} +O(\sm^2).
\end{equation}
As the probability of the bead is $b_{i_0}\rho_{\bv\wv}\sm+O(\sm
^2)$, the conditional probability we want is given by
%
%e72 ###
%
\begin{equation}\label{eq:probnontypical}
\prob{N\mbox{ successive $c$-edges}|\mbox{bead}}=\frac
{\rho_{\bv\wv}}{\rho_{\bv\wv_N}} \prod_{i=1}^{N} \frac
{c_{[i]}}{a_{[i-1]}} +O(\sm).
\end{equation}
%
% \begin{equation}
% A_{N+1}= [\begin{array}{c|cccc}\displaystyle
% o & \frac{1}{a_{[i_0]j}} & o & \cdots& o\\
% \hline
% o& x & \frac{1}{a_{[i_0+1]j}} & o & o \\
% \vdots& x& & \ddots& o\\
% o & & & & \frac{1}{a_{[i_0+N]j}}\\
% \sm\rho_{\bv\wv_N} +O(\sm^2)
% \end{array}]
% \end{equation}
Using Proposition \ref{lem:Qrank1}, one can rewrite $\rho_{\bv\wv
_N}/\rho_{\bv\wv}$ as the following telescopic product:
%
%e73 ###
%
\begin{equation}
\frac{\rho_{\bv\wv_N}}{\rho_{\bv\wv}}=\frac{U_{\bv}V_{\wv
_N}}{U_{\bv}V_{\wv}}=
\prod_{i=1}^{N}\frac{U_{\bv_i}V_{\wv_i}}{U_{\bv_i}V_{\wv_{i-1}}}.
\end{equation}
Plugging this into \eqref{eq:probnontypical}, one gets
%
%e74 ###
%
\begin{eqnarray}
\prob{N\mbox{ successive $c$-edges}|\mbox{bead}}&=&
\prod_{i=1}^N \frac{c_{[i]}U_{\bv_i}V_{\wv_i}}{a_{[i-1]}U_{\bv
_i}V_{\wv_{i-1}}}+O(\sm)\nonumber\\[-8pt]\\[-8pt]
&=&
\prod_{i=1}^N \frac{\ell(c_{[i]})}{\ell(a_{[i-1]})}+O(\sm),\nonumber
\end{eqnarray}
where $\ell(a_{[i-1]})$ and $\ell(c_{[i]})$ are respectively the
length of the dual edges with weight $a_{[i-1]}$ and $c_{[i]}$ given by
the mapping described in Lemma \ref{lem:mapping}.
In particular, in the limit $\sm\rightarrow0$, the probability that
the length $L$ of this succession of nontypical edges exceeds $p$
fundamental domains equals
%
%e75 ###
%
\begin{equation}
\prob{L\geq p}=
\Biggl(\prod_{i=1}^m \frac{\ell(c_{i})}{\ell(a_i)}\Biggr)^p.
\end{equation}
Thus, in the limit, $L$ has a geometric distribution.
\end{pf}

Pushing further the above computations of the lengths of the nontypical
sequences of edges in frozen columns, one can derive the following.

\begin{prop}
The limiting bead models on the different families of threads $\mathbf
{B}_j$ are perfectly correlated: The distance between beads on each
side of a frozen column converges in probability to zero.
\end{prop}
\begin{pf}
% Since the distribution of the lengths of sequences of non typical
%edges is close to the geometric distribution for a small $\sm$, the
%probability of
The details of the proof are omitted here, but by looking carefully at
the determinants in the proof given above, one can in fact see that,
for every $\sm$, the probability that a sequence of non typical edges
in a frozen columns exceeds, say, $\frac{1}{\sqrt{\sm}}$ is of order
$q^{1/\sqrt{\sm}}$, and thus decays very fast when $\sm$ goes to
zero. Thus, in the vertically rescaled graph, the distance between two
beads at the extremity of a sequence of nontypical edges is close to 0.
In the scaling limit, the distance between these beads converges in
probability to 0.
\end{pf}

\begin{figure}

\includegraphics{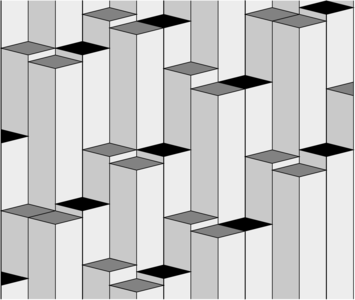}

\caption{A sketch of a typical tiling for a Gibbs
measure corresponding to a point inside a tentacle. The black beads
mark the transition between the two dominant type of edges. Next to
them, beads also appear on the other threads to compensate the defect
created in frozen columns.}
\label{fig:perles_gen}
\end{figure}

As a consequence, the picture of a typical dimer configuration for a
Gibbs measure corresponding to a point in a tentacle of the phase
diagram looks like the one in Figure \ref{fig:perles_gen}.

%f11 ###
%

\section*{Acknowledgments}

We are grateful to thank Richard Kenyon for the many fruitful
discussions and precious advice on this problem. We would like to
thank also Daniel Slutsky and Jon Warren for pointing me to reference
\cite{Barysh}, Scott Sheffield for showing me the random surface
interpretation of the bead model, as well as the referee for useful
comments.

%imsgetref loaded by dianan, Friday, October, 17:09:20 2008
%

\printaddresses

\end{document}